\documentclass[12pt,a4paper,twoside,reqno]{amsart}

\usepackage{babel}
\usepackage[utf8]{inputenc}
\usepackage[T1]{fontenc}
\usepackage{lmodern}
\usepackage{amsmath,amsfonts,amssymb,amsthm,mathtools}
\usepackage{geometry}
\usepackage{enumitem}
\usepackage{xcolor}
\usepackage{microtype}
\usepackage{url}
\usepackage{faktor}
\usepackage{mathrsfs}
\usepackage{eucal}
\usepackage{booktabs}
\usepackage{array}

\usepackage{listings}

\usepackage[
bookmarks=true,
breaklinks=true,
bookmarksnumbered=true,
colorlinks=true,
urlcolor=green,
anchorcolor=yellow,
citecolor=blue,
linkcolor=blue
]{hyperref}

\allowdisplaybreaks

\newcommand{\R}{\mathbb R}
\newcommand{\Z}{\mathbb Z}

\newcommand{\Pbar}{\overline P_n}

\newcommand{\Isom}{\operatorname{Isom}}

\newcommand{\rank}{\operatorname{rank}}
\newcommand{\Aut}{\operatorname{Aut}}
\newcommand{\Gi}[1]{\mathcal G_{#1}}
\newcommand{\Bi}[1]{\mathcal B_{#1}}
\newcommand{\XGi}[1]{X_{\Gi{#1}}}
\newcommand{\XBi}[1]{X_{\Bi{#1}}}
\newcommand{\Gammai}[2]{\Gamma^{(#1)}_{#2}}
\newcommand{\Ki}[2]{K^{(#1)}_{#2}}

\newcommand{\GammaB}[2]{\Gamma^{(\mathcal B_{#1})}_{#2}}
\newcommand{\KB}[2]{K^{(\mathcal B_{#1})}_{#2}}

\newtheorem{theorem}{Theorem}[section]
\newtheorem{proposition}[theorem]{Proposition}
\newtheorem{lemma}[theorem]{Lemma}

\newtheorem{theoremintro}{Theorem}

\theoremstyle{definition}

\newtheorem{example}[theorem]{Example}

\theoremstyle{remark}
\newtheorem{remark}[theorem]{Remark}

\newtheorem{remarkintro}[theoremintro]{Remark}

\title[{\tiny Crystallographic braid quotients and torus bundles over flat $3$-manifolds}]{Crystallographic quotients of Artin braid groups and torus bundles over flat $3$-manifolds}
\author[Oscar Ocampo]{Oscar Ocampo}
\address{Oscar Ocampo\\ Universidade Federal da Bahia, Departamento de Matem\'atica - IME, CEP:~40170-110 - Salvador, Brazil}
\email{oscaro@ufba.br}


\date{\today}

\begin{document}

\begin{abstract}
For $n\geq3$, let $\overline B_n=\faktor{B_n}{[P_n,P_n]}$ be the crystallographic quotient of the Artin braid group, with translation subgroup $\overline P_n=\faktor{P_n}{[P_n,P_n]}\cong\mathbb Z^N$, where $N=\binom n2$. We prove that every compact connected flat $3$-manifold $M$, considered up to affine equivalence, is the base of a geometric torus bundle whose total space is a compact flat manifold with fundamental group a Bieberbach subgroup of some $\overline B_n$. More precisely, for every $n\geq n_M$, where $n_M\in\{3,4,5\}$ is an explicit optimal threshold depending on $M$, we construct a Bieberbach subgroup $\Gamma_n(M)\leq\overline B_n$ of dimension $N$ and a complete normal subgroup $K_n(M)\cong\mathbb Z^{N-3}$ consisting of translations such that
$$
\faktor{\Gamma_n(M)}{K_n(M)}\cong\pi_1(M).
$$
We also determine which of the associated torus-bundle group extensions split: for the realizations constructed here, splitting occurs for exactly four of the ten affine-equivalence classes of compact flat $3$-manifolds. 
\end{abstract}

\maketitle

\section{Introduction}
\label{sec:intro}

The quotient of the Artin braid group by the commutator subgroup of the pure braid group, $\overline B_n=\faktor{B_n}{[P_n,P_n]}$, provides a natural bridge between braid groups and crystallographic groups. Indeed, there is a short exact sequence
$$
1 \longrightarrow \overline P_n \longrightarrow \overline B_n \overset{\overline\sigma}{\longrightarrow}
S_n \longrightarrow 1,
$$
where $\overline P_n:=\faktor{P_n}{[P_n,P_n]}$ is free abelian of rank $N=\binom n2=\frac{n(n-1)}2$. The group $\overline B_n$ is crystallographic of dimension $N$ with holonomy group $S_n$ \cite[Proposition~1]{GGO2017}, and the associated holonomy representation is the permutation action of $S_n$ on the set of unordered pairs of strings \cite[Proposition~12]{GGO2017}.

Bieberbach subgroups of $\overline B_n$ with prescribed finite holonomy have been constructed in several settings. If $H\leq S_n$ is a finite $2$-group, then the full preimage $\overline\sigma^{-1}(H)$ is a Bieberbach group \cite[Lemma~18]{GGO2017}. Bieberbach subgroups with cyclic holonomy of order $2^d$ were studied in \cite{OR2019}, and constructions with arbitrary finite abelian holonomy were obtained in \cite{Ocampo2021}. These results produce large families of compact flat manifolds naturally associated with quotients of Artin braid groups.

The question considered here is of a different nature. Prescribing the abstract holonomy group does not determine a flat manifold. One must also take into account the integral holonomy representation and the extension class of the corresponding Bieberbach group. We therefore ask whether the compact flat $3$-manifolds themselves, rather than merely their holonomy groups, can be recovered geometrically from the crystallographic braid quotients.

We use the notation of Wolf \cite[Theorems~3.5.5 and~3.5.9]{Wolf}. The six orientable $3$-dimensional Bieberbach groups are denoted by $\Gi{1},\Gi{2},\Gi{3},\Gi{4},\Gi{5},\Gi{6}$, and the four non-orientable ones by
$\Bi{1},\Bi{2},\Bi{3},\Bi{4}$. 
We write
$$
\XGi{i}=\R^3/\Gi{i}, \qquad \XBi{j}=\R^3/\Bi{j} 
$$
for the corresponding compact flat $3$-manifolds. Their holonomy groups are, respectively, $1, \mathbb Z_2, \mathbb Z_3, \mathbb Z_4, \mathbb Z_6, \mathbb Z_2^2$ in the orientable case, and $\mathbb Z_2, \mathbb Z_2, \mathbb Z_2^2, \mathbb Z_2^2$ in the non-orientable case. The manifold $\XGi{6}$ is the classical Hantzsche--Wendt manifold.

Our main result gives a braid-theoretic realization of all ten affine-equivalence classes of compact flat $3$-manifolds.

\begin{theoremintro}\label{thm:main}

Let $M$ be a compact connected flat $3$-manifold, considered up to affine equivalence. Then there exists an integer $n_M$ such that, for every $n\geq n_M$, there exist a Bieberbach subgroup $\Gamma_n(M)\leq\overline B_n$ of dimension $N=\binom n2$ and a complete normal subgroup $K_n(M)\cong\mathbb Z^{N-3}$ consisting of translations such that
$$
\faktor{\Gamma_n(M)}{K_n(M)}\cong\pi_1(M).
$$
Consequently, the compact flat manifold $X_{\Gamma_n(M)}=\R^N/\Gamma_n(M)$ admits a geometric fiber bundle
\begin{equation}\label{eq:main-fibration}
T^{N-3}\longrightarrow X_{\Gamma_n(M)}\longrightarrow X_{\faktor{\Gamma_n(M)}{K_n(M)}},
\end{equation}
whose fibers are totally geodesic flat tori, and whose base is affinely equivalent to $M$.
\end{theoremintro}

\begin{remarkintro}\label{rmk:mainthm}
The values in Theorem~\ref{thm:main}, supplied by the constructions in this paper, may be chosen as follows, and they are optimal:
$$
\begin{array}{c|cccccccccc}
M
&
\XGi{1}
&
\XGi{2}
&
\XGi{3}
&
\XGi{4}
&
\XGi{5}
&
\XGi{6}
&
\XBi{1}
&
\XBi{2}
&
\XBi{3}
&
\XBi{4}
\\ \hline
n_M
&
3&4&3&4&5&4&3&3&4&4.
\end{array}
$$

Indeed, the condition $N=\binom n2\geq3$ forces $n\geq3$. For the remaining lower bounds, recall that the holonomy group of $\pi_1(M)\cong\faktor{\Gamma_n(M)}{K_n(M)}$ is a quotient of the holonomy group $F\leq S_n$ of $\Gamma_n(M)$; see Lemma~\ref{lem:holonomy-quotient}. If $\pi_1(M)$ has holonomy group $\mathbb Z_4$ or $\mathbb Z_2^2$, then $4$ divides $|F|$, which is impossible for $F\leq S_3$; hence $n\geq4$. If $\pi_1(M)$ has holonomy group $\mathbb Z_6$, then $F$ admits an abelian quotient of order $6$, so $6$ divides $|F|$ and $F^{\mathrm{ab}}$ maps onto $\mathbb Z_6$. The subgroups of $S_4$ of order divisible by $6$ are the conjugates of $S_3$, the alternating group $A_4$ and $S_4$ itself, whose abelianizations are $\mathbb Z_2$, $\mathbb Z_3$ and $\mathbb Z_2$, respectively; hence $n\geq5$. These observations give all the required lower bounds except for $\Gi{2}$. If $n=3$, then $N=3$, so $K_3(M)$ is trivial and $\Gamma_3(M)\cong\pi_1(M)$; in particular, the holonomy group of $\pi_1(M)$ is a subgroup of $S_3$, acting on the three unordered pairs. Every non-trivial involution of $S_3$ acts on these pairs as a transposition, and hence with determinant $-1$. Thus any $3$-dimensional Bieberbach subgroup of $\overline B_3$ with holonomy group $\mathbb Z_2$ is non-orientable, whereas $\Gi{2}$ is orientable. Therefore $n_{\XGi{2}}\geq4$ as well.
\end{remarkintro}

The theorem is stronger than a realization of the possible abstract holonomy groups. At the algebraic level, the basic mechanism is the following. Starting with the translation lattice $\Lambda$ of a Bieberbach subgroup $\Gamma\leq\overline B_n$, we construct a sublattice $K\leq\Lambda$ of rank $N-3$, normal in $\Gamma$, such that $\faktor{\Lambda}{K}\cong\mathbb Z^3$. The induced holonomy action on this rank-$3$ quotient determines the integral holonomy representation of the base, while the relations satisfied by the chosen holonomy lifts determine the corresponding Bieberbach extension. Thus both the integral representation and the extension data enter essentially into the construction.

This distinction is especially visible for the two non-orientable groups $\Bi{3}$ and $\Bi{4}$. Both are obtained from the Klein four subgroup 
$$
V^-=\langle(1,2),(3,4)\rangle\leq S_4
$$
and from a common finite-index translation lattice. Their $3$-dimensional quotients have the same integral holonomy representation, but the chosen holonomy lifts satisfy different mixed relations, yielding different extension classes. In this way the braid-group construction distinguishes $\Bi{3}$ from $\Bi{4}$ even though neither their abstract holonomy group nor their integral holonomy representation distinguishes them.

The embedding of the holonomy group in the symmetric group is itself geometrically relevant. For instance, the Hantzsche--Wendt case $\Gi{6}$ uses
$$
V^+ = \langle (1,4)(2,3),(1,3)(2,4) \rangle \leq S_4, 
$$
whose three non-trivial elements are double transpositions. In contrast, $V^-$ contains two transpositions. The two Klein four subgroups are therefore not conjugate in $S_4$, and they act differently on the six unordered pairs: $V^+$ has three orbits of length $2$, whereas $V^-$ fixes $\{1,2\}$ and $\{3,4\}$ and permutes the remaining four pairs freely.

The ten realizations require several complementary integral constructions, including stabilization, finite-index corrections of translation lattices, cyclotomic module quotients, and suitable choices of holonomy lifts. 

Beyond the existence of the torus fibrations, we determine exactly which of the associated torus-bundle group extensions split. The splitting classification is established through integral holonomy-module obstructions and explicit relations among the chosen holonomy lifts.

For the explicit realizations constructed below, we write $\Gammai{i}{n}=\Gamma_n(\XGi{i})$ for $2\leq i\leq6$ and $\GammaB{j}{n}=\Gamma_n(\XBi{j})$ for $1\leq j\leq4$, referring to the particular choices of groups in Theorem~1. We use analogous notation for the corresponding kernels. 
A particularly illustrative phenomenon is that, for $n\geq4$, 
$$
\GammaB{2}{n}=\Gammai{2}{n}
$$
as subgroups of $\overline B_n$, as follows directly from the constructions in Sections~\ref{subsec:G2} and~\ref{subsec:B12}. Nevertheless, the corresponding torus-bundle extension over $\XBi{2}$ splits, whereas the one over $\XGi{2}$ does not. Thus the torus-bundle structure retains information about the chosen complete normal subgroup that is not determined by the total-space Bieberbach group alone.

The trivial-holonomy case deserves a brief comment. Since $\Gi{1}\cong\mathbb Z^3$ and $\XGi{1}=T^3$, for every $n\geq3$ one may take $\Gamma_n(\XGi{1}) = \overline P_n \cong \mathbb Z^N$ and choose any rank-$(N-3)$ direct summand $K\leq\overline P_n$ with $\faktor{\overline P_n}{K} \cong \mathbb Z^3$. This gives the product bundle 
$$
T^{N-3} \longrightarrow T^N \longrightarrow T^3. 
$$
For the other cases with $n_M=3$, the fiber in \eqref{eq:main-fibration} is $T^0$ at the minimal value $n=3$, so the corresponding $3$-dimensional Bieberbach group already occurs as a subgroup of $\overline B_3$.

The explicit nature of these realizations also makes them a convenient source of examples for questions in topological fixed point theory. For instance, compact flat manifolds are infra-nilmanifolds, and the averaging formula of Kim, Lee and Lee~\cite[Theorem~3.5]{KLL2005} computes the Nielsen number of any self-map of such a manifold from the linear parts of an affine homotopy representative and the holonomy; the torus bundles constructed here provide natural test cases for Nielsen and Reidemeister numbers of fiber-preserving maps. Similarly, the regular double coverings in Proposition~\ref{prop:G4-covering} and Remark~\ref{rem:B1B2-covering} give explicit free involutions on higher-dimensional flat manifolds, suggesting the study of their Borsuk--Ulam properties in the spirit of the work of Bauval, Gon\c{c}alves and Hayat~\cite{BGH} in dimension three. We do not pursue these questions here.

We conclude the Introduction with an outline of the paper. Section~\ref{sec:preliminaries} recalls the required facts about crystallographic groups and geometric fibrations and records a faithful affine realization of $\overline B_n$. Section~\ref{sec:cyclic} treats the five non-trivial orientable cases: the cyclic holonomy cases $\Gi{2}, \ldots, \Gi{5}$ and the Hantzsche--Wendt case $\Gi{6}$. Section~\ref{sec:nonorientable} treats the four non-orientable cases. Section~\ref{sec:realization-splitting} completes the proof of Theorem~\ref{thm:main} and classifies the associated torus-bundle extensions as split or non-split.  

\subsection*{Acknowledgments}

O.~O.~was partially supported by the National Council for Scientific and Technological Development (CNPq, Brazil) through a \textit{Bolsa de Produtividade} grant No.~305422/2022--7.

\section{Preliminaries}\label{sec:preliminaries}

This section collects the background and establishes the general framework used throughout the paper. We first recall the relation between complete normal subgroups of Bieberbach groups and geometric fibrations. We then review the crystallographic braid quotient $\overline B_n$ and its affine realization, and finally fix Wolf's notation for the ten compact flat $3$-manifolds.

\subsection{Crystallographic groups and geometric fibrations}\label{subsec:bieb}

Recall that $\Isom(\R^m)=O(m)\ltimes\R^m$. A crystallographic group of dimension $m$ is a discrete cocompact subgroup of $\Isom(\R^m)$, and a torsion-free crystallographic group is called a Bieberbach group; see \cite{Charlap,Szczepanski,Wolf}.

Let $\Gamma$ be a Bieberbach group of dimension $m$. Its translation subgroup $\Lambda$ is free abelian of rank $m$ and is normal and maximal abelian in $\Gamma$. Moreover, there is a short exact sequence $1\longrightarrow\Lambda\longrightarrow\Gamma\longrightarrow F\longrightarrow1$, where $F$ is the holonomy group. Conjugation induces a faithful integral representation $\rho\colon F\longrightarrow\Aut(\Lambda)\cong GL(m,\Z)$.  The quotient $X_\Gamma=\R^m/\Gamma$ is a compact flat Riemannian manifold with fundamental group $\Gamma$. 

For later use, recall that a subgroup $K\leq\Lambda\cong\mathbb Z^m$ is called \emph{primitive} if $\faktor{\Lambda}{K}$ is torsion-free; equivalently, $K$ is a direct summand of $\Lambda$. An element $v\in\Lambda$ is called primitive if the cyclic subgroup $\langle v\rangle$ is primitive.

Following Ratcliffe and Tschantz~\cite{RT2010}, we use the term \emph{space group} for a crystallographic group. We use their theory of complete normal subgroups of space groups. We do not need the full definition of completeness here. The two facts required below are that a normal subgroup $\mathscr{C}$ of a space group $\Gamma$ is complete if and only if $\faktor{\Gamma}{\mathscr{C}}$ is again a space group \cite[Theorem~5]{RT2010}, and that, when $\faktor{\Gamma}{\mathscr{C}}$ is torsion-free, the corresponding geometric fibration is a geometric fiber bundle over the flat manifold associated with $\faktor{\Gamma}{\mathscr{C}}$ \cite[Theorem~13]{RT2010}.

\begin{proposition}\label{prop:fibration-principle}
Let $\Gamma$ be a Bieberbach group of dimension $N$ with translation lattice $\Lambda$, and let $K\leq\Lambda$ be a normal free abelian subgroup of rank $N-3$. Suppose that $\faktor{\Gamma}{K}$ is a $3$-dimensional Bieberbach group. Then $K$ is complete and the compact flat manifold $X_\Gamma$ admits a geometric fiber bundle
$$
T^{N-3}\longrightarrow X_\Gamma\longrightarrow X_{\faktor{\Gamma}{K}},
$$
with totally geodesic fibers.
\end{proposition}

\begin{proof}
Since $\faktor{\Gamma}{K}$ is a Bieberbach group, it is in particular a space group. Hence $K$ is complete by \cite[Theorem~5]{RT2010}. Since $\faktor{\Gamma}{K}$ is torsion-free, \cite[Theorem~13]{RT2010} then implies that $X_\Gamma$ admits a geometric fiber bundle over the flat manifold associated with $\faktor{\Gamma}{K}$, with fiber $V/K$, where $V=\operatorname{Span}_{\mathbb R}(K)$. 

Since $K\leq\Lambda$ consists of translations and has rank $N-3$, one has $\dim V=N-3$ and $V/K\cong T^{N-3}$. The base is precisely $X_{\faktor{\Gamma}{K}}$, and the fibers, being the images of affine cosets of $V$, are totally geodesic.
\end{proof}

The next lemma compares the holonomy groups of $\Gamma$ and of $\faktor{\Gamma}{K}$.

\begin{lemma}\label{lem:holonomy-quotient}
Let $\Gamma$ be a Bieberbach group with translation lattice $\Lambda$ and holonomy group $F$, and let $K\leq\Lambda$ be a normal subgroup of $\Gamma$ such that $\Pi=\faktor{\Gamma}{K}$ is a Bieberbach group. Let $L$ and $F_\Pi$ denote the translation lattice and the holonomy group of $\Pi$, and let $p\colon\Gamma\longrightarrow\Pi$ be the quotient map. Then $p(\Lambda)=\faktor{\Lambda}{K}\leq L$, and $p$ induces a surjective homomorphism $F\longrightarrow F_\Pi$ whose kernel is isomorphic to $\faktor{L}{p(\Lambda)}$. In particular, if $|F|=|F_\Pi|$, then $p(\Lambda)=L$, the restriction $p|_\Lambda\colon\Lambda\longrightarrow L$ is surjective with kernel $K$, and $p$ induces an isomorphism $F\cong F_\Pi$.
\end{lemma}

\begin{proof}
Set $A=p(\Lambda)$. Since $K\leq\Lambda$, we have $\faktor{\Pi}{A}\cong\faktor{\Gamma}{\Lambda}=F$. In particular, $A$ is an abelian normal subgroup of finite index in $\Pi$.

We first show that $A$ consists entirely of translations. Since $A$ has finite index in $\Pi$, the subgroup $A\cap L$ has finite index in $L$. Hence its translation vectors span $\mathbb R^d$, where $d=\rank L$. Let $\gamma=(R,v)\in A$. Since $A$ is abelian, $\gamma$ commutes with every translation $(I,w)\in A\cap L$. Thus
$$
\gamma(I,w)\gamma^{-1}=(I,Rw)=(I,w).
$$
It follows that $Rw=w$ for every such $w$, and therefore $R=I$. Hence $A\leq L$.

Consequently, the inclusions $A\leq L\leq\Pi$ give an exact sequence
$$
1\longrightarrow\faktor{L}{A}\longrightarrow\faktor{\Pi}{A}\longrightarrow\faktor{\Pi}{L}\longrightarrow1.
$$
Since $\faktor{\Pi}{A}\cong F$ and $\faktor{\Pi}{L}=F_\Pi$, this yields
$$
1\longrightarrow\faktor{L}{p(\Lambda)}\longrightarrow F\longrightarrow F_\Pi\longrightarrow1.
$$
In particular, if $|F|=|F_\Pi|$, then $L=p(\Lambda)$ and the induced homomorphism $F\longrightarrow F_\Pi$ is an isomorphism. The remaining assertions follow from $\ker(p)=K\leq\Lambda$.
\end{proof}

\begin{example}
Let $H=\langle a,b\mid aba^{-1}=b^{-1}\rangle$ be the fundamental group of the Klein bottle, and consider the $4$-dimensional Bieberbach group $\Gamma=H\times\langle c\rangle\times\langle d\rangle$. Its translation lattice is $\Lambda=\langle a^2,b,c,d\rangle\cong\mathbb Z^4$, and its holonomy group is $F\cong\mathbb Z_2$.

For $K_1=\langle c\rangle$, the quotient $\Pi_1=\Gamma/K_1\cong H\times\mathbb Z$ is a $3$-dimensional Bieberbach group with holonomy $F_{\Pi_1}\cong\mathbb Z_2$. Thus the induced homomorphism $F\longrightarrow F_{\Pi_1}$ is an isomorphism.

On the other hand, for $K_2=\langle b\rangle$, the quotient $\Pi_2=\Gamma/K_2\cong\mathbb Z^3$ has trivial holonomy. If $p_2\colon\Gamma\longrightarrow\Pi_2$ denotes the quotient map and $L_2$ is the translation lattice of $\Pi_2$, then $p_2(\Lambda)$ has index $2$ in $L_2$. Consequently, $L_2/p_2(\Lambda)\cong\mathbb Z_2$, which is precisely the kernel of $F\longrightarrow F_{\Pi_2}$.
\end{example}

\subsection{The crystallographic braid quotient}\label{subsec:braid}

Let $B_n$ be the Artin braid group on $n$ strings, with standard generators $\sigma_1,\ldots,\sigma_{n-1}$. Let $\sigma\colon B_n\longrightarrow S_n$ be the natural homomorphism defined by $\sigma(\sigma_i)=(i,i+1)$, and let $P_n=\ker(\sigma)$. We follow the convention of \cite{GGO2017} and read permutations from left to right. 

For $1\leq i<j\leq n$, let
\begin{equation}
A_{i,j}=\sigma_{j-1}\cdots\sigma_{i+1}\sigma_i^2 \sigma_{i+1}^{-1}\cdots\sigma_{j-1}^{-1}.
\label{eq:Aij}
\end{equation}
As in \cite{GGO2017}, whenever necessary we use the convention $A_{j,i}=A_{i,j}$. The quotient $\overline P_n=\faktor{P_n}{[P_n,P_n]}$ is free abelian of rank $N=\binom n2$, with basis $\{A_{i,j}\mid1\leq i<j\leq n\}$. By abuse of notation, the image of $A_{i,j}$ in $\overline P_n$ will again be denoted by $A_{i,j}$. 

Set $\overline B_n=\faktor{B_n}{[P_n,P_n]}$. Then
$$
1\longrightarrow\overline P_n\longrightarrow\overline B_n \overset{\overline\sigma}{\longrightarrow}S_n\longrightarrow1
$$
is a crystallographic extension \cite[Proposition~1]{GGO2017}. If $\gamma\in\overline B_n$ and $\pi=\overline\sigma(\gamma^{-1})$, then
\begin{equation}
\gamma A_{i,j}\gamma^{-1}=A_{\pi(i),\pi(j)}.
\label{eq:pair-action}
\end{equation}
Thus the holonomy representation is the permutation representation associated with the natural action of $S_n$ on unordered pairs \cite[Proposition~12]{GGO2017}.

We next record a concrete faithful affine realization of $\overline B_n$, which will also be used in the explicit computations below. Let $E=\mathbb R^N$, $N=\binom n2$, with orthonormal basis
$$
\{e_{i,j}\mid 1\leq i<j\leq n\}.
$$
As above, whenever necessary we use the convention $e_{j,i}=e_{i,j}$. We write an affine isometry of $E$ as $(P,b)$, acting by $(P,b)(x)=Px+b$.  Thus
$$
(P,b)(Q,c)=(PQ,Pc+b).
$$

For $1\leq i\leq n-1$, let $P_i$ be the permutation matrix induced by the adjacent transposition $(i,i+1)$ on the set of unordered pairs, and define
$$
S_i = \left(P_i,\frac12e_{i,i+1}\right) \in\Isom(E).
$$

\begin{proposition}\label{prop:affine-model}
For $n\geq3$, the assignment $\sigma_i\longmapsto S_i$ induces an injective homomorphism $\Phi_n\colon\overline B_n\longrightarrow\Isom(E)$. Moreover, $\Phi_n(A_{i,j})=(I,e_{i,j})$, so the restriction of $\Phi_n$ to $\overline P_n$ identifies $\overline P_n$ with the lattice of translations generated by the vectors $e_{i,j}$.

If $\gamma\in\overline B_n$ and $\pi=\overline\sigma(\gamma^{-1})$, then the linear part of $\Phi_n(\gamma)$ is the permutation matrix induced by $e_{i,j}\longmapsto e_{\pi(i),\pi(j)}$. 
\end{proposition}

\begin{proof}
We first verify that the assignment $\sigma_i\mapsto S_i$ defines a homomorphism on $B_n$. If $|i-j|\geq2$, then $P_iP_j=P_jP_i$, and $P_i$ fixes $e_{j,j+1}$ while $P_j$ fixes $e_{i,i+1}$. Therefore
$$
S_iS_j = S_jS_i.
$$
For adjacent indices, the linear parts satisfy $P_iP_{i+1}P_i = P_{i+1}P_iP_{i+1}$. Using the affine multiplication formula, the translation part of both $S_iS_{i+1}S_i$ and $S_{i+1}S_iS_{i+1}$ is $\frac12 \left( e_{i,i+1}+e_{i,i+2}+e_{i+1,i+2} \right)$. Hence 
$$
S_iS_{i+1}S_i=S_{i+1}S_iS_{i+1}.
$$
Thus the elements $S_i$ satisfy the defining relations of the Artin braid group and determine a homomorphism $\widehat\Phi_n\colon B_n\longrightarrow\Isom(E)$. 

Since $P_i$ fixes $e_{i,i+1}$, one has  $S_i^2=(I,e_{i,i+1})$. More generally, for every affine isometry $(P,b)$ and every translation $t_v=(I,v)$, 
$$
(P,b)t_v(P,b)^{-1}=t_{Pv}.
$$
Using the definition of $A_{i,j}$ in \eqref{eq:Aij}, it follows that  $\widehat\Phi_n(A_{i,j})=(I,e_{i,j})$. In particular, $P_n$ is mapped into the abelian group of translations, so 
$$
[P_n,P_n]\leq\ker(\widehat\Phi_n).
$$
Consequently, $\widehat\Phi_n$ factors through a homomorphism
$$
\Phi_n\colon\overline B_n\longrightarrow\Isom(E).
$$

The restriction of $\Phi_n$ to $\overline P_n=\faktor{P_n}{[P_n,P_n]}$ sends the basis $\{A_{i,j}\}$ onto the standard lattice basis $\{e_{i,j}\}$. Hence it is injective on $\overline P_n$. 

With our convention for multiplication in $S_n$, the linear part of $\Phi_n(\gamma)$, for $\gamma\in\overline B_n$, is the permutation matrix associated with $\pi=\overline\sigma(\gamma^{-1})$, which is precisely the action appearing in \eqref{eq:pair-action}. 

Finally, suppose that $\Phi_n(\gamma)=1$. Its linear part is trivial. For $n\geq3$, the natural action of $S_n$ on unordered pairs is faithful, so $\overline\sigma(\gamma)=1$. Hence $\gamma\in\overline P_n$. Since $\Phi_n$ is injective on $\overline P_n$, we conclude that $\gamma=1$. Therefore $\Phi_n$ is injective.
\end{proof}

Henceforth we identify $\overline P_n$ with the integer translation lattice
$$
\bigoplus_{1\leq i<j\leq n}\mathbb Z e_{i,j} \subset E.
$$

We repeatedly use the fact that $\overline B_n$ has no non-trivial finite-order element of even order \cite[Theorem~2]{GGO2017}. In particular:

\begin{proposition}[{\cite[Lemma~18]{GGO2017}}]\label{prop:2groups}
If $H\leq S_n$ is a finite $2$-group, then $\overline\sigma^{-1}(H)\leq \overline B_n$ is a Bieberbach group of dimension $N$ with holonomy group $H$.
\end{proposition}

\begin{remark}\label{rem:stabilization}
Let $3\leq m\leq n$. The standard inclusion $B_m\hookrightarrow B_n$ induces an embedding $\overline B_m\hookrightarrow\overline B_n$; see \cite[Theorem~7]{Ocampo2021} and compare \cite[Theorem~3(a)]{GGO2017}. This also follows from Proposition~\ref{prop:affine-model}. We identify $\overline B_m$ with its image, so that $\overline B_m\cap\overline P_n=\overline P_m$. Set
$$
E^{(m)}_n=\langle A_{i,j}\mid 1\leq i<j\leq n,\ j\geq m+1\rangle.
$$
Then $\overline P_n=\overline P_m\oplus E^{(m)}_n$, and $E^{(m)}_n$ is normalized by $\overline B_m$ by \eqref{eq:pair-action}. Consequently, for every subgroup $\Gamma_m\leq\overline B_m$, we have $\Gamma_m\cap E^{(m)}_n=\{1\}$ and hence
$$
\Gamma_mE^{(m)}_n=E^{(m)}_n\rtimes\Gamma_m,\qquad \faktor{\Gamma_mE^{(m)}_n}{E^{(m)}_n}\cong\Gamma_m.
$$
   All the stabilized groups in Sections~\ref{sec:cyclic} and~\ref{sec:nonorientable} are of this form, with $E_{3,n}=E^{\mathcal B}_n=E^{(3)}_n$, $E_{2,n}=E_{4,n}=E^-_n=E^{(4)}_n$ and $E_{5,n}=E^{(5)}_n$; the $\Gi{6}$ construction also uses $E^{(4)}_n$.
\end{remark}

\subsection{Wolf's notation for flat 3-manifolds}\label{subsec:wolf}

For later comparison with the quotient groups constructed below, we record presentations of the orientable Bieberbach groups with non-trivial cyclic holonomy. These presentations are equivalent to those in \cite[Theorem~3.5.5]{Wolf}. 

The group $\Gi{2}$ has translation subgroup $\langle z,u,v\rangle\cong\Z^3$ and a lift $t$ of the non-trivial element of $\mathbb Z_2$ satisfying
\begin{equation}
t^2=z,\qquad tut^{-1}=u^{-1},\qquad tvt^{-1}=v^{-1}. 
\label{eq:G2-presentation}
\end{equation}

The group $\Gi{3}$ has holonomy $\mathbb Z_3$ and admits a presentation with commuting translations $z,u,v$ and a holonomy lift $t$ satisfying
$$
t^3=z,\qquad tut^{-1}=v,\qquad tvt^{-1}=u^{-1}v^{-1}. 
$$

The group $\Gi{4}$ has holonomy $\mathbb Z_4$ and admits a presentation with commuting translations $z,u,v$ and a holonomy lift $t$ satisfying 
\begin{equation}
t^4=z,\qquad tut^{-1}=v,\qquad tvt^{-1}=u^{-1}.
\label{eq:G4-presentation}
\end{equation}

The group $\Gi{5}$ has holonomy $\mathbb Z_6$ and admits a presentation with commuting translations $z,u,v$ and a holonomy lift $t$ satisfying
\begin{equation}
t^6=z,\qquad tut^{-1}=v,\qquad tvt^{-1}=u^{-1}v. 
\label{eq:G5-presentation}
\end{equation}

The remaining orientable case, $\Gi{6}$, has holonomy $\mathbb Z_2^2$, and $\XGi{6}$ is the classical Hantzsche--Wendt manifold. Unlike the cyclic cases above, its identification in Section~\ref{sec:G6} relies on the uniqueness of the orientable $3$-dimensional Bieberbach group with this holonomy, rather than on an explicit presentation. Presentations of the four non-orientable groups $\Bi{1},\ldots,\Bi{4}$ will be used explicitly in Section~\ref{sec:nonorientable}; see also \cite[List~(A)]{HJKL2002}.

\section{The orientable flat $3$-manifolds}\label{sec:cyclic}

We now construct braid-theoretic realizations of the five non-trivial orientable Bieberbach groups $\Gi{2},\ldots,\Gi{6}$. The four cyclic holonomy cases are treated first, followed by the Hantzsche--Wendt case, whose holonomy is the Klein four group. Although the constructions require different integral ingredients, they all fit into the quotient-to-fibration framework of Section~\ref{sec:preliminaries}.

Throughout this section, for $2\leq i\leq6$, whenever $\Gammai{i}{n}$ is defined, we write
$$
X_{i,n} = X_{\Gammai{i}{n}} = \faktor{\mathbb R^N}{\Gammai{i}{n}}.
$$

\subsection{The case $\Gi{2}$: holonomy $\mathbb Z_2$}\label{subsec:G2}

We begin with the simplest non-trivial orientable case, for which the full preimage of a transposition already provides the required total-space Bieberbach group. Let $n\geq4$ and let $H_2=\langle(1,2)\rangle\leq S_n$. By Proposition~\ref{prop:2groups}, 
$$
\Gammai{2}{n} := \widetilde H_{2,n} := \overline\sigma^{-1}(H_2)
$$
is a Bieberbach group of dimension $N$ and holonomy $\mathbb Z_2$. 

For $n=4$, \cite[Example~8(b)]{OR2019} shows that 
\begin{equation}
\faktor{\widetilde H_{2,4}}{K_{2,4}} \cong \Gi{2},
\label{eq:G2-seed}
\end{equation}
where $K_{2,4} = \langle A_{1,3}A_{2,3}, A_{1,4}A_{2,4}, A_{3,4} \rangle \cong \mathbb Z^3$. Here the quotient is identified with the $3$-dimensional Bieberbach group $\Gi{2}$ in Wolf's classification.

For $n\geq4$, set
\begin{equation}
E_{2,n}=\langle A_{i,j}\mid 1\leq i<j\leq n,\ j\geq5\rangle,
\label{eq:E2n}
\end{equation} 
where $E_{2,4}=\{1\}$. This is an $H_2$-invariant free abelian subgroup of rank $N-6$. Define
$$
\Ki{2}{n}=K_{2,4}\oplus E_{2,n}.
$$

\begin{theorem}\label{thm:G2}
For every $n\geq4$, the subgroup $\Ki{2}{n}$ is normal in $\Gammai{2}{n}$, has rank $N-3$, and $\faktor{\Gammai{2}{n}}{\Ki{2}{n}}\cong\Gi{2}$.  Consequently, 
$$
T^{N-3}\longrightarrow X_{2,n}\longrightarrow\XGi{2}
$$
is a geometric torus bundle with totally geodesic fibers.
\end{theorem}

\begin{proof}
The subgroup $E_{2,n}$ is invariant under the transposition $(1,2)$: for every $j\geq5$ it interchanges $A_{1,j}$ and $A_{2,j}$ and preserves the set of all other generators occurring in \eqref{eq:E2n}. The three generators of $K_{2,4}$ are fixed by the holonomy action and hence are central in $\Gammai{2}{4}$. Hence $\Ki{2}{n}$ is normal. 

Quotienting by $E_{2,n}$ removes all translation coordinates involving strings $5,\ldots,n$ and reduces the construction to the $n=4$ case. Quotienting further by $K_{2,4}$ gives \eqref{eq:G2-seed}. Since 
$$
\rank\Ki{2}{n}=3+(N-6)=N-3,
$$
the group-theoretic assertion follows. The fibration statement is Proposition~\ref{prop:fibration-principle}. 
\end{proof}

\subsection{The case $\Gi{3}$: holonomy $\mathbb Z_3$}\label{subsec:G3}

The full preimage of a subgroup of order $3$ in $S_n$ need not be torsion-free, so the construction differs from the preceding case. We start from the Bieberbach subgroup exhibited in \cite[Remarks~26(c)]{GGO2017}. 

In $\overline B_3$, put
$$
x=A_{1,2}\sigma_1^{-1}\sigma_2, \qquad s=A_{1,2}A_{1,3}A_{2,3}. 
$$
Then $x^3=s$. Let $\Lambda_{3,3}=\langle s,A_{1,3}^3,A_{2,3}^3\rangle$. The subgroup $\Gammai{3}{3}=\langle\Lambda_{3,3},x\rangle$ is a $3$-dimensional Bieberbach group with holonomy $\mathbb Z_3$ \cite[Remarks~26(c)]{GGO2017}. Its holonomy representation is orientation preserving; hence, by Wolf's classification, 
\begin{equation}
\Gammai{3}{3}\cong\Gi{3}. 
\label{eq:G3-identification}
\end{equation}

For $n\geq3$, regard the same elements as elements of $\overline B_n$ via the embedding of Remark~\ref{rem:stabilization}, and set $E_{3,n}=\langle A_{i,j}\mid 1\leq i<j\leq n,\ j\geq4\rangle$. The subgroup $E_{3,n}$ is invariant under the $3$-cycle induced by $x$. Define 
$$
\Lambda_{3,n}=\Lambda_{3,3}\oplus E_{3,n}, \qquad \Gammai{3}{n}=\langle\Lambda_{3,n},x\rangle.
$$

\begin{theorem}\label{thm:G3}
For every $n\geq3$, the group $\Gammai{3}{n}$ is a Bieberbach subgroup of $\overline B_n$ of dimension $N$ and holonomy $\mathbb Z_3$. Moreover, $\faktor{\Gammai{3}{n}}{E_{3,n}}\cong\Gi{3}$.  Consequently,
$$
T^{N-3} \longrightarrow X_{3,n} \longrightarrow \XGi{3}
$$
is a geometric torus bundle with totally geodesic fibers. 
\end{theorem}

\begin{proof}
The lattice $\Lambda_{3,n}$ has full rank $N$ in $\overline P_n$, and the conjugation formula \eqref{eq:pair-action} shows that it is normalized by $x$. The quotient group $\faktor{\Gammai{3}{n}}{E_{3,n}}$ is exactly the $n=3$ group $\Gammai{3}{3}$ of \eqref{eq:G3-identification}. 

If $g\in\Gammai{3}{n}$ has finite order, its image in $\Gammai{3}{3}$ has finite order and is therefore trivial. Hence $g\in E_{3,n}$, which is free abelian, and so $g=1$. Thus $\Gammai{3}{n}$ is torsion-free. Since it contains the full-rank lattice $\Lambda_{3,n}$ with finite quotient $\mathbb Z_3$, it is crystallographic and hence a Bieberbach group. The last assertion follows from Proposition~\ref{prop:fibration-principle}.
\end{proof}

\subsection{The case $\Gi{4}$: holonomy $\mathbb Z_4$}\label{subsec:G4}

The case of holonomy $\mathbb Z_4$ requires an integral modification of the natural translation lattice. We first carry out the construction for $n=4$.

Let
$$
\theta=\sigma_3\sigma_2\sigma_1\in\overline B_4, \qquad \mu=\overline\sigma(\theta), \qquad H_4=\langle\mu\rangle\leq S_4.
$$
Then $\mu$ is a $4$-cycle and $H_4\cong\mathbb Z_4$. The full preimage $\overline\sigma^{-1}(H_4)$ is a Bieberbach group by Proposition~\ref{prop:2groups}. Proposition~9 of \cite{OR2019}, specialized to $d=2$ and $m=1$, gives  
\begin{equation}
\theta^4=\prod_{1\leq i<j\leq4}A_{i,j}.
\label{eq:G4-fourth}
\end{equation}

Under conjugation by $\theta$, the six generators $A_{i,j}$ split into one orbit of length $4$ and one orbit of length $2$. Choose $x_0=A_{1,2}$,  $x_{r+1}=\theta x_r\theta^{-1}$ for $r\in\Z/4\Z$, and let $y_0=A_{1,3}$,  $y_1=\theta y_0\theta^{-1}$. Thus
$$
\{x_0,x_1,x_2,x_3\}=\{A_{1,2},A_{2,3},A_{3,4},A_{1,4}\}, \qquad \{y_0,y_1\}=\{A_{1,3},A_{2,4}\}.
$$
Put $s=y_0y_1$,  $d=y_0y_1^{-1}$. Then $\theta s\theta^{-1}=s$ and $\theta d\theta^{-1}=d^{-1}$. The lattice
$$
\Lambda_{4,4}=\langle x_0,x_1,x_2,x_3,s,d\rangle
$$
has index $2$ in $\overline P_4$. Indeed, replacing the basis $(y_0,y_1)$ by $(s,d)$ has determinant $-2$. Moreover, \eqref{eq:G4-fourth} becomes 
\begin{equation}
\theta^4=x_0x_1x_2x_3s\in\Lambda_{4,4}.
\label{eq:G4-fourth-new}
\end{equation}
Define
$$
\Gammai{4}{4}=\langle\Lambda_{4,4},\theta\rangle.
$$
It is an index-$2$ subgroup of the Bieberbach group $\overline\sigma^{-1}(H_4)$ and is therefore a Bieberbach group of dimension $6$ and holonomy $\mathbb Z_4$.

Now set
\begin{equation}
K_{4,4}=\langle x_0x_2,\ x_1x_3,\ d\rangle.
\label{eq:K44}
\end{equation}
The subgroup $K_{4,4}$ is normal in $\Gammai{4}{4}$ because
$$
\theta(x_0x_2)\theta^{-1}=x_1x_3, \qquad \theta(x_1x_3)\theta^{-1}=x_0x_2, \qquad \theta d\theta^{-1}=d^{-1}.
$$
It is primitive in $\Lambda_{4,4}$: in additive notation its generators are $x_0+x_2$, $x_1+x_3$, and $d$, and the quotient is free abelian of rank $3$. Write $u=\overline{x}_0$, $v=\overline{x}_1$, $z=\overline{s}$ for the corresponding basis of $\faktor{\Lambda_{4,4}}{K_{4,4}}$. Then 
\begin{equation}
\overline\theta u\overline\theta^{-1}=v, \qquad \overline\theta v\overline\theta^{-1}=u^{-1}, \qquad \overline\theta z\overline\theta^{-1}=z,
\label{eq:G4-action-quot}
\end{equation}
and \eqref{eq:G4-fourth-new} gives
\begin{equation}
\overline\theta^4=z.
\label{eq:G4-primitive}
\end{equation}

\begin{theorem}
\label{thm:G4}
For every $n\geq4$, there exists a Bieberbach subgroup $\Gammai{4}{n}\leq\overline B_n$ of dimension $N$ and holonomy $\mathbb Z_4$, together with a normal subgroup $\Ki{4}{n}\cong\Z^{N-3}$, such that $\faktor{\Gammai{4}{n}}{\Ki{4}{n}}\cong\Gi{4}$. Consequently, 
$$
T^{N-3} \longrightarrow X_{4,n} \longrightarrow \XGi{4} 
$$
is a geometric torus bundle with totally geodesic fibers.
\end{theorem}

\begin{proof}
For $n=4$, put $Q_4=\faktor{\Gammai{4}{4}}{K_{4,4}}$.  Since $K_{4,4}$ is contained in the translation lattice, $Q_4$ fits into the short exact sequence 
$$
1\longrightarrow \faktor{\Lambda_{4,4}}{K_{4,4}} \longrightarrow Q_4 \longrightarrow \mathbb Z_4 \longrightarrow1, 
$$
where $\faktor{\Lambda_{4,4}}{K_{4,4}} \cong \langle u,v,z\rangle \cong\Z^3$. Relations~\eqref{eq:G4-action-quot}--\eqref{eq:G4-primitive} give 
$$
\overline\theta^4=z,\qquad \overline\theta u\overline\theta^{-1}=v,\qquad \overline\theta v\overline\theta^{-1}=u^{-1},\qquad \overline\theta z\overline\theta^{-1}=z.
$$
Together with commutativity of $u,v,z$, these are precisely the defining relations for the Bieberbach group $\Gi{4}$ in Wolf's classification; compare \eqref{eq:G4-presentation}. Hence $Q_4\cong\Gi{4}$. 

For $n\geq 5$, let $E_{4,n} = \langle A_{i,j}\mid 1\leq i<j\leq n,\ j\geq5\rangle$. This subgroup is invariant under the $4$-cycle. Set 
$$
\Lambda_{4,n} = \Lambda_{4,4}\oplus E_{4,n}, \qquad \Gammai{4}{n} = \langle\Lambda_{4,n},\theta\rangle, 
$$
and $\Ki{4}{n} = K_{4,4}\oplus E_{4,n}$. The lattice $\Lambda_{4,n}$ has index $2$ in $\overline P_n$, and $\Gammai{4}{n}$ has index $2$ in the Bieberbach group $\overline\sigma^{-1}(H_4)$. Hence $\Gammai{4}{n}$ is a Bieberbach group. Moreover, 
$$
\rank\Ki{4}{n} = 3+(N-6) = N-3,
$$
and quotienting by $E_{4,n}$ reduces the construction to the $n=4$ case. Therefore
$$
\faktor{\Gammai{4}{n}}{\Ki{4}{n}}\cong\Gi{4}.
$$
Proposition~\ref{prop:fibration-principle} gives the geometric fibration.
\end{proof}

\begin{remark}\label{rem:G4-integral}
The replacement $(y_0,y_1)\mapsto(s,d)$ is the essential integral correction in this case. To see why it is needed, consider the quotient of $\langle\overline P_4,\theta\rangle$ by the normal subgroup $\langle x_0x_2,\ x_1x_3,\ y_0y_1^{-1}\rangle$. Writing $\bar y$ for the common image of $y_0$ and $y_1$, equation~\eqref{eq:G4-fourth} gives $\theta^4=\bar y^2$. Since $\theta^2$ fixes $\bar y$, the element $\theta^2\bar y^{-1}$ satisfies $(\theta^2\bar y^{-1})^2=1$. Moreover, this element is non-trivial, since its image under the holonomy projection is $\mu^2\neq1$. Thus the resulting quotient contains an element of order $2$ and is not torsion-free.

By contrast, passing to the index-$2$ lattice $\Lambda_{4,4}$ makes $s=y_0y_1$ a primitive basis element. Its image $z$ is then primitive in the rank-$3$ quotient $\faktor{\Lambda_{4,4}}{K_{4,4}}$, and the relation becomes $\theta^4=z$, as required for $\Gi{4}$.
\end{remark}

\begin{proposition}\label{prop:G4-covering}
For every $n\geq4$, let $\widetilde H_{4,n} = \overline\sigma^{-1}(H_4)$.  Then $\Gammai{4}{n}\triangleleft\widetilde H_{4,n}$, and $\faktor{\widetilde H_{4,n}}{\Gammai{4}{n}} \cong\mathbb Z_2$.  Consequently, there is a regular double Riemannian covering
$$
X_{4,n} \longrightarrow X_{\widetilde H_{4,n}} 
$$
with deck transformation group $\mathbb Z_2$.
\end{proposition}

\begin{proof}
Since both groups map surjectively onto $H_4$ and $\theta\in\Gammai{4}{n}$, we have
$$
\widetilde H_{4,n} = \overline P_n\Gammai{4}{n}.
$$
We first prove that $\Gammai{4}{n}$ is normal in $\widetilde H_{4,n}$. Since $\overline P_n$ is abelian, it centralizes the translation lattice $\Lambda_{4,n}$, so it is enough to show that $\overline P_n$ normalizes the coset $\Lambda_{4,n}\theta$. 

For $p\in\overline P_n$,
$$
p\theta p^{-1}\theta^{-1} = p\bigl(\theta p^{-1}\theta^{-1}\bigr). 
$$
The subgroup generated by these elements is contained in $\Lambda_{4,n}$. Indeed, on the orbit $\{x_0,x_1,x_2,x_3\}$ the corresponding differences are generated by the elements $x_i x_{i+1}^{-1}$, all of which belong to $\Lambda_{4,n}$. On the orbit $\{y_0,y_1\}$ they are generated by $y_0y_1^{-1}=d$, and every translation involving an index at least $5$ already belongs to $E_{4,n}\leq\Lambda_{4,n}$. Hence $p\theta p^{-1}\in\Lambda_{4,n}\theta \subseteq \Gammai{4}{n}$, and therefore $\Gammai{4}{n}\triangleleft\widetilde H_{4,n}$. 

Now the second isomorphism theorem gives
$$
\faktor{\widetilde H_{4,n}}{\Gammai{4}{n}} \cong \faktor{\overline P_n}{\overline P_n\cap\Gammai{4}{n}} = \faktor{\overline P_n}{\Lambda_{4,n}} \cong \mathbb Z_2.
$$
The inclusion of Bieberbach groups acting on the same Euclidean space induces the stated regular double Riemannian covering.
\end{proof}

\subsection{The case $\Gi{5}$: holonomy $\mathbb Z_6$}\label{subsec:G5}

We now specialize the finite abelian holonomy construction of \cite[Theorem~1(a)(ii)]{Ocampo2021}. In the notation used there, we specialize the parameters denoted by $t$ and $v$ in \cite{Ocampo2021} to $1$, and take $p_1=3$ and $r_1=r_2=1$. Then $G=\mathbb Z_3\times\mathbb Z_2\cong\mathbb Z_6$ embeds in $S_5$ by a disjoint $3$-cycle and $2$-cycle, and \cite{Ocampo2021} gives a Bieberbach subgroup of $\overline B_5$ of dimension $10$ with this holonomy. For the reader's convenience, we give a self-contained proof of the properties of this group that we need; see Lemma~\ref{lem:G5-seed}.

Let $a=A_{1,2}\delta_{0,3}$ and $b=\sigma_4$, where $\delta_{0,3}=\sigma_2\sigma_1^{-1}$ is the element of \cite[Eq.~(8)]{Ocampo2021} (see also \cite[Eq.~(14)]{GGO2017}). It has order $3$ in $\overline B_n$ \cite[Lemma~28]{GGO2017}, and $\overline\sigma(\delta_{0,3})$ is a $3$-cycle on $\{1,2,3\}$. The elements $a$ and $b$ involve disjoint blocks of strings, and hence commute. Put
$$
s=a^3=A_{1,2}A_{1,3}A_{2,3}, \qquad t=b^2=A_{4,5}.
$$
Indeed, $a^3=A_{1,2}\,\bigl(\delta_{0,3}A_{1,2}\delta_{0,3}^{-1}\bigr)\bigl(\delta_{0,3}^2A_{1,2}\delta_{0,3}^{-2}\bigr)\,\delta_{0,3}^3$, where $\delta_{0,3}^3=1$ and, by \eqref{eq:pair-action}, the three displayed conjugates of $A_{1,2}$ are $A_{1,2}$, $A_{1,3}$ and $A_{2,3}$ in some order. Consider the lattice
\begin{equation}
\mathcal L_{5,5} = \left\langle s,\ A_{1,3}^3,\ A_{2,3}^3; A_{1,4},A_{1,5},A_{2,4},A_{2,5},A_{3,4},A_{3,5};  t \right\rangle,
\label{eq:G5-lattice}
\end{equation}
where the displayed generators form a basis; in the present case, this is the basis $C_1\cup C_2\cup C_3$ of \cite[Eq.~(14)]{Ocampo2021}. Let $\theta=ab$. Its image $\overline\sigma(\theta)\in S_5$, the product of a $3$-cycle on $\{1,2,3\}$ and the transposition $(4,5)$, has order $6$. Moreover, since $a$ and $b$ commute,
\begin{equation}
\theta^6=a^6b^6=s^2t^3.
\label{eq:G5-sixth}
\end{equation}

\begin{lemma}\label{lem:G5-seed}
The group $\Gammai{5}{5}=\langle\mathcal L_{5,5},\theta\rangle\leq\overline B_5$ is a Bieberbach group of dimension $10$ with translation lattice $\mathcal L_{5,5}$ and holonomy group $\langle\overline\sigma(\theta)\rangle\cong\mathbb Z_6$. It coincides with the Bieberbach group of \cite[Theorem~1(a)(ii)]{Ocampo2021} for the above choice of parameters.
\end{lemma}

\begin{proof}
By \eqref{eq:pair-action}, conjugation by $\theta$ acts on $\overline P_5$ through a permutation of $\{1,\ldots,5\}$ that permutes $1,2,3$ cyclically and interchanges $4$ and $5$. Hence it fixes $s$ and $t$, permutes the six basis elements $A_{i,j}$ with $1\leq i\leq3<j\leq5$, and permutes the set $\{A_{1,2}^3,A_{1,3}^3,A_{2,3}^3\}$, which is contained in $\mathcal L_{5,5}$ because $A_{1,2}^3=s^3A_{1,3}^{-3}A_{2,3}^{-3}$. Thus $\mathcal L_{5,5}$ is normalized by $\theta$. Since $\theta^6\in\mathcal L_{5,5}$ by \eqref{eq:G5-sixth}, while $\overline\sigma(\theta^k)\neq1$ for $1\leq k\leq5$, we obtain
$$
\Gammai{5}{5}=\bigsqcup_{k=0}^{5}\mathcal L_{5,5}\,\theta^k \qquad\text{and}\qquad \Gammai{5}{5}\cap\overline P_5=\mathcal L_{5,5}.
$$
As $\mathcal L_{5,5}$ has rank $10$, the group $\Gammai{5}{5}$ is crystallographic of dimension $10$, with translation lattice $\mathcal L_{5,5}$ and holonomy group $\langle\overline\sigma(\theta)\rangle\cong\mathbb Z_6$.

It remains to show that $\Gammai{5}{5}$ is torsion-free; in order to do this we adapt the argument of \cite[Remarks~26(c)]{GGO2017}. If $w\in\Gammai{5}{5}$ is a non-trivial element of finite order, then $\langle w\rangle\cap\mathcal L_{5,5}=\{1\}$ because $\mathcal L_{5,5}$ is torsion-free, so $\langle w\rangle$ maps injectively into $\mathbb Z_6$; moreover, since $\Gamma_5^{(5)}\leq\overline B_5$, the order of $w$ is odd by \cite[Theorem~2]{GGO2017}. Hence it suffices to show that $\Gammai{5}{5}$ has no element of order $3$. Such an element projects to an element of order $3$ of $\langle\overline\sigma(\theta)\rangle$, so it lies in $\mathcal L_{5,5}\theta^2\cup\mathcal L_{5,5}\theta^4$. Since $\theta^2=a^2t$ and $\theta^4=a\,st^2$, and since $s$ and $t$ commute with $a$, it has the form $w=\lambda a^j$ with $\lambda\in\mathcal L_{5,5}$ and $j\in\{1,2\}$. Write $\overline P_5$ additively, and let $g$ denote the automorphism of $\overline P_5$ induced by conjugation by $a$. Since $a^3=s$ is a translation, $g^3=1$, and therefore
$$
w^3=\lambda+g^{j}\lambda+g^{2j}\lambda+js=\lambda+g\lambda+g^2\lambda+js.
$$
Let $f\colon\overline P_5\longrightarrow\mathbb Z$ be the sum of the coordinates of $A_{1,2}$, $A_{1,3}$ and $A_{2,3}$ with respect to the basis $\{A_{i,j}\}$. Since $g$ permutes $\{A_{1,2},A_{1,3},A_{2,3}\}$ and permutes the remaining basis elements, $f\circ g=f$. Moreover, $f(s)=f(A_{1,3}^3)=f(A_{2,3}^3)=3$ and $f$ vanishes on the other basis elements of $\mathcal L_{5,5}$, so $f(\mathcal L_{5,5})\subseteq3\mathbb Z$. If $w^3=1$, then $3f(\lambda)+3j=0$, that is, $f(\lambda)=-j$, which is impossible because $f(\lambda)\in3\mathbb Z$ and $j\in\{1,2\}$. Hence $\Gammai{5}{5}$ is torsion-free.

Finally, the group of \cite[Theorem~1(a)(ii)]{Ocampo2021} is generated by $a$, $b$, the elements $A_{r,s}^3$ with $1\leq r<s\leq3$, and the elements $A_{i,j}$ with $j\geq4$. These elements generate the same subgroup as $\mathcal L_{5,5}$ and $\theta$, since $\theta=ab$, $s=a^3$, $a=\theta^4s^{-1}t^{-2}$ and $b=a^{-1}\theta$.
\end{proof}

We now construct an invariant subgroup of rank $7$ in the translation lattice. First, in $L_0=\langle s,p,q\rangle,$ with $p=A_{1,3}^3$ and $q=A_{2,3}^3$, define
$$
\ell\colon L_0\longrightarrow\Z, \qquad \ell(s)=\ell(p)=\ell(q)=1. 
$$
Equivalently, $\ell$ is induced by one third of the sum of the exponent coordinates of $A_{1,2},A_{1,3},A_{2,3}$. Since the $\mathbb Z_3$-action permutes these three coordinates, $\ell$ is holonomy invariant. Hence
$$
K_0=\ker(\ell)=\langle ps^{-1},qs^{-1}\rangle
$$
is invariant and has rank $2$, while $\faktor{L_0}{K_0}\cong\Z$ with primitive generator $\overline s$.

Second, consider the rank-$6$ cross-term lattice 
$$
W=\langle A_{1,4},A_{1,5},A_{2,4},A_{2,5},A_{3,4},A_{3,5}\rangle.
$$
The generator $\theta$ acts transitively on these six basis elements. Choose $w_0=A_{1,4}$ and write $w_j=\theta^jw_0\theta^{-j}$, for $0\leq j\leq5$. Thus, writing the translation lattice additively, $W\cong\Z[\Z_6]$ as a $\Z[\Z_6]$-module. Let
\begin{equation}
q_6\colon W\longrightarrow \faktor{\Z[x]}{(x^2-x+1)}
\label{eq:G5-cyclotomic-map}
\end{equation}
be the $\Z[\Z_6]$-module epimorphism determined by $q_6(w_j)=x^j$, where the generator of $\Z_6$ acts on the target by multiplication by $x$. 

Put $K_W=\ker(q_6)$. Since $x^2-x+1=\Phi_6(x)$ is monic and divides $x^6-1$, the quotient $\faktor{W}{K_W}$ is free abelian of rank $2$. In the quotient $\faktor{\Gamma_5^{(5)}}{K_W}$, let $u$ and $v$ denote the images of $w_0$ and $w_1$, respectively, and let $\overline\theta$ denote the image of $\theta$. Then
\begin{equation}
\overline\theta u\overline\theta^{-1}=v, \qquad \overline\theta v\overline\theta^{-1}=u^{-1}v. 
\label{eq:G5-plane-action}
\end{equation}

Finally, in the fixed rank-$2$ lattice $\langle s,t\rangle$, impose the relation
\begin{equation}
st=1.
\label{eq:G5-axis-relation}
\end{equation}
The element $st$ is primitive in $\langle s,t\rangle$. Thus put $K_{5,5}=\langle K_0,K_W,st\rangle$. The quotient of the rank-$4$ lattice $\langle s,p,q,t\rangle$ by $\langle ps^{-1},qs^{-1},st\rangle$ is infinite cyclic, generated by $\overline s$. Therefore $K_{5,5}$ is primitive of rank $2+4+1=7$ in the rank-$10$ lattice \eqref{eq:G5-lattice}. Set $z=\overline s^{-1}$. Then \eqref{eq:G5-sixth} and \eqref{eq:G5-axis-relation} give 
\begin{equation}
\overline\theta^6=\overline{s}^2\overline{t}^3 = \overline{s}^{-1}=z.
\label{eq:G5-primitive-sixth}
\end{equation}

\begin{theorem}\label{thm:G5} 
For every $n\geq5$, there exists a Bieberbach subgroup $\Gammai{5}{n}\leq\overline B_n$ of dimension $N$ and holonomy $\mathbb Z_6$, together with a normal subgroup $\Ki{5}{n}\cong\Z^{N-3}$, such that $\faktor{\Gammai{5}{n}}{\Ki{5}{n}}\cong\Gi{5}$. Consequently,
$$
T^{N-3} \longrightarrow X_{5,n} \longrightarrow \XGi{5}
$$
is a geometric torus bundle with totally geodesic fibers.
\end{theorem}

\begin{proof}
For $n=5$, the group $\Gammai{5}{5}=\langle \mathcal L_{5,5},\theta\rangle$ is a Bieberbach group with translation lattice $\mathcal L_{5,5}$ and holonomy group $\mathbb Z_6$, by Lemma~\ref{lem:G5-seed}.
The subgroups $K_0$ and $K_W$ are invariant by construction, and $s$ and $t$ are fixed by the holonomy. Hence $K_{5,5}\triangleleft\Gammai{5}{5}$. 

Put $Q_5=\faktor{\Gammai{5}{5}}{K_{5,5}}$. The quotient fits into
$$
1\longrightarrow\Z^3 \longrightarrow Q_5 \longrightarrow\mathbb Z_6 \longrightarrow1, 
$$
where the translation subgroup has basis $z,u,v$. Relations~\eqref{eq:G5-plane-action} and~\eqref{eq:G5-primitive-sixth} give 
$$
\overline\theta^6=z,\qquad \overline\theta u\overline\theta^{-1}=v,\qquad \overline\theta v\overline\theta^{-1}=u^{-1}v,
$$
and $z,u,v$ commute. These are precisely the defining relations of $\Gi{5}$ in Wolf's classification; compare \eqref{eq:G5-presentation}. Therefore $Q_5\cong\Gi{5}$. 

For $n\geq 6$, let
$$
E_{5,n} = \langle A_{i,j}\mid1\leq i<j\leq n,\ j\geq6\rangle. 
$$
The $\mathbb Z_6$-action fixes the indices $6,\ldots,n$, so $E_{5,n}$ is invariant. Define
$$
\mathcal L_{5,n}=\mathcal L_{5,5}\oplus E_{5,n}, \qquad \Gammai{5}{n}=\langle \mathcal L_{5,n},\theta\rangle,
$$
and $\Ki{5}{n}=K_{5,5}\oplus E_{5,n}$. The quotient $\faktor{\Gammai{5}{n}}{E_{5,n}}$ is $\Gammai{5}{5}$. Hence any torsion element of $\Gammai{5}{n}$ would lie in $E_{5,n}$ and must be trivial. Thus $\Gammai{5}{n}$ is torsion-free; since $\mathcal L_{5,n}$ is a full-rank translation lattice with finite cyclic quotient, it is a Bieberbach group. Finally,
$$
\rank\Ki{5}{n} = 7+(N-10) = N-3,
$$
and the quotient by $\Ki{5}{n}$ is unchanged. Hence $\faktor{\Gammai{5}{n}}{\Ki{5}{n}}\cong\Gi{5}$.  Proposition~\ref{prop:fibration-principle} completes the proof.
\end{proof}

\subsection{The case $\Gi{6}$ and the Klein four group}\label{sec:G6}

We now treat the unique orientable flat $3$-manifold with non-cyclic holonomy. Let $n\geq4$ and put
$$
g_1=(1,4)(2,3), \qquad g_2=(1,3)(2,4), \qquad G=\langle g_1,g_2\rangle\cong\mathbb Z_2^2. 
$$
Define
\begin{align}
d_1&=A_{2,3}A_{1,4}^{-1}, & s_1&=A_{1,2}A_{3,4},\nonumber\\
d_2&=A_{2,4}A_{1,3}^{-1}, & s_2&=A_{1,3}A_{2,4},\nonumber\\
d_3&=A_{3,4}A_{1,2}^{-1}, & s_3&=A_{1,4}A_{2,3}.\nonumber
\end{align}
Let $\Lambda_{6,n}$ be the sublattice of $\overline P_n$ generated by 
$$
\mathcal B_{6,n}=
\left( d_1,d_2,d_3,s_1,s_2,s_3;  (A_{1,j},A_{4,j},A_{2,j},A_{3,j})_{5\leq j\leq n};  (A_{i,j})_{5\leq i<j\leq n} \right).
$$
On each of the three rank-$2$ blocks involving the first four strings, the corresponding sum--difference generators span a sublattice of index $2$. Hence 
\begin{equation}
[\overline P_n:\Lambda_{6,n}]=8, \qquad \faktor{\overline P_n}{\Lambda_{6,n}}\cong\mathbb Z_2^3. 
\label{eq:G6-index}
\end{equation}
The lattice $\Lambda_{6,n}$ is $G$-invariant. On $\langle d_1,d_2,d_3\rangle$, the action is
\begin{equation}
\rho(g_1)=\operatorname{diag}(1,-1,-1), \qquad \rho(g_2)=\operatorname{diag}(-1,1,-1), 
\label{eq:G6-3action}
\end{equation}
while $s_1,s_2,s_3$ are fixed by $G$.

Choose lifts
$$
\alpha=A_{2,3}\sigma_3\sigma_2\sigma_1\sigma_2\sigma_3\sigma_2, \qquad \beta=A_{2,4}\sigma_2\sigma_1\sigma_2^2\sigma_3\sigma_2.
$$
Then $\overline\sigma(\alpha)=g_1$, $\overline\sigma(\beta)=g_2$, and a direct computation using the affine model of Proposition~\ref{prop:affine-model} gives
\begin{align}
\alpha^2&=d_1s_1s_2s_3^2, \label{eq:G6-alpha2}\\
\beta^2&=d_2s_1s_2^2s_3, \label{eq:G6-beta2}\\
(\alpha\beta)^2&=d_3s_1^2s_2^3s_3^3. \label{eq:G6-alphabeta2}
\end{align}
Define
$$
\Gammai{6}{n}=\langle\Lambda_{6,n},\alpha,\beta\rangle.
$$

\begin{proposition}\label{prop:Gamma6-bieb}
For every $n\geq4$, the group $\Gammai{6}{n}$ is a Bieberbach group of dimension $N$ with translation lattice $\Lambda_{6,n}$ and holonomy group $G\cong\mathbb Z_2^2$. 
\end{proposition}

\begin{proof}
The $G$-invariance of $\Lambda_{6,n}$ and relations \eqref{eq:G6-alpha2}--\eqref{eq:G6-alphabeta2} show that $\Lambda_{6,n}\triangleleft\Gammai{6}{n}$ and that the quotient by $\Lambda_{6,n}$ is generated by two commuting involutions. Since the induced map onto $G$ is surjective, it is an isomorphism. Hence the translation subgroup is exactly $\Lambda_{6,n}$.

The full preimage $\overline\sigma^{-1}(G)$ is a Bieberbach group by Proposition~\ref{prop:2groups}, so its subgroup $\Gammai{6}{n}$ is discrete and torsion-free. Since $\Lambda_{6,n}$ has full rank in $\overline P_n$, the quotient $\mathbb R^N/\Gammai{6}{n}$ is compact. Thus $\Gammai{6}{n}$ is a Bieberbach group with the claimed dimension and holonomy.
\end{proof}

Define the $G$-invariant sublattices
$$
\Lambda_{6,n}^{\mathrm{base}}=\langle d_1,d_2,d_3\rangle
$$
and
$$
\Ki{6}{n} = \langle s_1,s_2,s_3, A_{i,j}\mid 1\leq i<j\leq n,\ j\geq5\rangle.
$$
Then
$$
\Lambda_{6,n}=\Lambda_{6,n}^{\mathrm{base}}\oplus\Ki{6}{n}, \qquad \rank\Ki{6}{n}=N-3=\frac{(n-3)(n+2)}2.
$$
Modulo $\Ki{6}{n}$, relations \eqref{eq:G6-alpha2}--\eqref{eq:G6-alphabeta2} become
\begin{equation}
\overline\alpha^2=\overline d_1, \qquad \overline\beta^2=\overline d_2, \qquad (\overline\alpha\,\overline\beta)^2=\overline d_3.
\label{eq:G6-square-quotient}
\end{equation}

\begin{theorem}\label{thm:G6}
For every $n\geq4$, $\faktor{\Gammai{6}{n}}{\Ki{6}{n}} \cong \Gi{6}$.  Consequently, 
$$
T^{N-3} \longrightarrow X_{6,n} \longrightarrow \XGi{6} 
$$
is a geometric torus bundle with totally geodesic fibers.
\end{theorem}

\begin{proof}
The quotient fits into
$$
1\longrightarrow \faktor{\Lambda_{6,n}}{\Ki{6}{n}} \longrightarrow \faktor{\Gammai{6}{n}}{\Ki{6}{n}} \longrightarrow G\longrightarrow1,
$$
where the translation subgroup is free abelian of rank $3$ with basis $\overline d_1,\overline d_2,\overline d_3$. The holonomy action is given by \eqref{eq:G6-3action}, so it is faithful and orientation preserving.

It remains to check torsion-freeness. If an element projects to $g_1$, it has the form
$$
x=\overline d_1^a\overline d_2^b\overline d_3^c\overline\alpha,
$$
and \eqref{eq:G6-3action}--\eqref{eq:G6-square-quotient} give
$$
x^2=\overline d_1^{2a+1}\neq1.
$$
Similarly, an element projecting to $g_2$ has square $\overline d_2^{2b+1}$, and an element projecting to $g_1g_2$ has square $\overline d_3^{2c+1}$. Hence the quotient is torsion-free.

We have therefore obtained an orientable $3$-dimensional Bieberbach group with holonomy $\mathbb Z_2^2$. By \cite[Theorem~3.5.5]{Wolf}, there is a unique such group up to affine equivalence, namely $\Gi{6}$. Hence $\faktor{\Gammai{6}{n}}{\Ki{6}{n}} \cong \Gi{6}$. The corresponding flat manifold $\XGi{6}$ is the classical Hantzsche--Wendt manifold. Proposition~\ref{prop:fibration-principle} gives the geometric fiber bundle. 
\end{proof}

\subsection{The index-$8$ covering in the $\Gi{6}$ case}\label{subsec:G6-covering} 

The lattice correction in \eqref{eq:G6-index} also gives a regular finite covering. This parallels the double covering obtained in Proposition~\ref{prop:G4-covering} from the index-$2$ lattice correction in the $\Gi{4}$ case, but here the lattice has index $8$ and the deck transformation group is $\mathbb Z_2^3$. 

\begin{proposition}\label{prop:G6-covering} 
Let $G=\langle g_1,g_2\rangle\cong\mathbb Z_2^2$ be the Klein four subgroup of $S_n$ introduced in Section~\ref{sec:G6}, and put $\widetilde G_n=\overline\sigma^{-1}(G)$. Then $\Gammai{6}{n}$ is normal in $\widetilde G_n$ and $\faktor{\widetilde G_n}{\Gammai{6}{n}}\cong\mathbb Z_2^3$. Consequently, there is a regular $8$-sheeted Riemannian covering $X_{6,n}\longrightarrow X_{\widetilde G_n}$ with deck transformation group $\mathbb Z_2^3$.
\end{proposition}

\begin{proof}
Both groups map surjectively onto $G$, hence $\widetilde G_n=\Pbar\Gammai{6}{n}$. We first show that $\Gammai{6}{n}$ is normal in $\widetilde G_n$. Since $\Pbar$ is abelian, it suffices to prove that $\Pbar$ normalizes $\Gammai{6}{n}$.

Let $p\in\Pbar$ and $\gamma\in\Gammai{6}{n}$, and put $g=\overline\sigma(\gamma)\in G$. Writing $g\cdot p$ for the holonomy action on $\Pbar$, we have
$$
p\gamma p^{-1}\gamma^{-1}=p(g\cdot p)^{-1}.
$$
We claim that this element belongs to $\Lambda_{6,n}$. Indeed, $G$ permutes the six basis elements of $\Pbar$ involving the first four strings in three pairs, and the corresponding differences belong to $\langle d_1,d_2,d_3\rangle\leq\Lambda_{6,n}$. Moreover, every basis element $A_{i,j}$ with $j\geq5$, together with its image under $G$, belongs to $\Lambda_{6,n}$. Since $\Pbar$ is abelian, the map $p\mapsto p(g\cdot p)^{-1}$ is a homomorphism. Therefore $p(g\cdot p)^{-1}\in\Lambda_{6,n}$ for every $p\in\Pbar$. It follows that $p\gamma p^{-1}\in\Gammai{6}{n}$, proving that $\Gammai{6}{n}\trianglelefteq\widetilde G_n$.

By Proposition~\ref{prop:Gamma6-bieb}, we have $\Gammai{6}{n}\cap\Pbar=\Lambda_{6,n}$. Hence the second isomorphism theorem and \eqref{eq:G6-index} give
$$
\faktor{\widetilde G_n}{\Gammai{6}{n}}\cong\faktor{\Pbar}{\Lambda_{6,n}}\cong\mathbb Z_2^3.
$$
In particular, $[\widetilde G_n:\Gammai{6}{n}]=8$.

Finally, by Proposition~\ref{prop:affine-model}, both groups act by Euclidean isometries on the same space $\R^N$. Since they are Bieberbach groups, their actions are free, properly discontinuous and cocompact. Thus the inclusion $\Gammai{6}{n}\leq\widetilde G_n$ induces a Riemannian covering
$$
X_{6,n}=\R^N/\Gammai{6}{n}\longrightarrow\R^N/\widetilde G_n=X_{\widetilde G_n}
$$
of degree $8$. Since $\Gammai{6}{n}$ is normal in $\widetilde G_n$, the covering is regular, with deck transformation group $\faktor{\widetilde G_n}{\Gammai{6}{n}}\cong\mathbb Z_2^3$.
\end{proof}

\section{The non-orientable flat $3$-manifolds}\label{sec:nonorientable}

We now turn to the four non-orientable Bieberbach groups $\Bi{1},\Bi{2},\Bi{3},\Bi{4}$ in Wolf's classification \cite[Theorem~3.5.9]{Wolf}. No new abstract holonomy groups are required: the essential point is to realize the appropriate integral holonomy representations and, in the last two cases, the appropriate Bieberbach extension classes. For comparison with the groups constructed below, we use the presentations recorded in \cite[List~(A)]{HJKL2002}, which are equivalent to those in Wolf's classification.

\subsection{The groups $\Bi{1}$ and $\Bi{2}$: holonomy $\mathbb Z_2$}\label{subsec:B12}

We begin with the two non-orientable manifolds having holonomy $\mathbb Z_2$. Both arise from the same transposition subgroup, but the $\Bi{1}$ construction requires an index-$2$ modification of the translation lattice.  Let $n\geq3$ and consider $H_{\mathcal B} = \langle(1,2)\rangle \leq S_n$. Set
$$
\GammaB{2}{n} = \overline\sigma^{-1}(H_{\mathcal B}).
$$
By Proposition~\ref{prop:2groups}, $\GammaB{2}{n}$ is a Bieberbach group of dimension $N=\binom n2$ and holonomy $\mathbb Z_2$.

For $n=3$, let $\varepsilon=\sigma_1\in\overline B_3$. Then $\varepsilon^2=A_{1,2}$, and the holonomy action is
$$
\varepsilon A_{1,2}\varepsilon^{-1}=A_{1,2}, \qquad \varepsilon A_{1,3}\varepsilon^{-1}=A_{2,3}, \qquad \varepsilon A_{2,3}\varepsilon^{-1}=A_{1,3}.
$$
By \cite[Theorem~25(b)]{GGO2017},
\begin{equation}
\GammaB{2}{3}\cong\Bi{2}.
\label{eq:B2-identification}
\end{equation}

We now obtain the other non-orientable Bieberbach group with holonomy $\mathbb Z_2$ from an index-$2$ modification of the translation lattice. Put $s=A_{1,3}A_{2,3}$, $d=A_{1,3}A_{2,3}^{-1}$, and let 
$$
\Lambda_{\mathcal B_1,3} = \langle A_{1,2},s,d\rangle. 
$$
The change of basis from $(A_{1,3},A_{2,3})$ to $(s,d)$ has matrix
$$
\begin{pmatrix}
1&1\\
1&-1
\end{pmatrix},
$$
whose determinant is $-2$. Hence $[\overline P_3:\Lambda_{\mathcal B_1,3}]=2$. Moreover,
\begin{equation}
\varepsilon s\varepsilon^{-1}=s, \qquad \varepsilon d\varepsilon^{-1}=d^{-1}, \qquad \varepsilon^2=A_{1,2}. 
\label{eq:B1-action}
\end{equation}
Thus $\Lambda_{\mathcal B_1,3}$ is $\varepsilon$-invariant. 

Define
$$
\GammaB{1}{3} = \langle\Lambda_{\mathcal B_1,3},\varepsilon\rangle.
$$
The group $\GammaB{1}{3}$ has index $2$ in $\GammaB{2}{3}$, and therefore is a $3$-dimensional Bieberbach group. The translation subgroup of $\GammaB{1}{3}$ is $\langle A_{1,2},s,d\rangle\cong\mathbb Z^3$, and the quotient by this subgroup is generated by the class of $\varepsilon$. Relations \eqref{eq:B1-action}, together with the commutativity of the translations, coincide with the extension presentation of $\Bi{1}$ in \cite[List~(A)]{HJKL2002}. Hence 
\begin{equation}
\GammaB{1}{3}\cong\Bi{1}.
\label{eq:B1-identification}
\end{equation}

For the stabilization, define $E^{\mathcal B}_n = \langle A_{i,j} \mid 1\leq i<j\leq n,\ j\geq4 \rangle$. We take $E^{\mathcal B}_3=\{1\}$. The subgroup $E^{\mathcal B}_n$ is invariant under $(1,2)$ and has rank $N-3$. For $n\geq3$, put $\Lambda_{\mathcal B_1,n} = \Lambda_{\mathcal B_1,3}\oplus E^{\mathcal B}_n$ and define
$$
\GammaB{1}{n} = \langle \Lambda_{\mathcal B_1,n},\varepsilon \rangle. 
$$
The lattice $\Lambda_{\mathcal B_1,n}$ has index $2$ in $\overline P_n$, so $\GammaB{1}{n}$ has index $2$ in $\GammaB{2}{n}$. In particular, it is a Bieberbach group of dimension $N$.

\begin{theorem}\label{thm:B12} 
For every $n\geq3$, the groups $\GammaB{1}{n}$ and $\GammaB{2}{n}$ are Bieberbach subgroups of $\overline B_n$ of dimension $N$ and holonomy $\mathbb Z_2$. Moreover, $\KB{1}{n} = \KB{2}{n} = E^{\mathcal B}_n \cong \mathbb Z^{N-3}$ is normal in both groups, and 
$$
\faktor{\GammaB{1}{n}}{\KB{1}{n}} \cong\Bi{1}, \qquad \faktor{\GammaB{2}{n}}{\KB{2}{n}} \cong\Bi{2}. 
$$
Consequently, there are geometric torus bundles $T^{N-3} \longrightarrow X_{\GammaB{1}{n}} \longrightarrow \XBi{1}$ and $T^{N-3} \longrightarrow X_{\GammaB{2}{n}} \longrightarrow \XBi{2}$, with totally geodesic fibers.
\end{theorem}

\begin{proof}
Normality of $E^{\mathcal B}_n$ follows from its invariance under the transposition $(1,2)$. Notice that $\rank E^{\mathcal B}_n=N-3$. Quotienting by $E^{\mathcal B}_n$ removes precisely the translation coordinates involving strings $4,\ldots,n$. Hence
$$
\faktor{\GammaB{1}{n}}{E^{\mathcal B}_n} \cong \GammaB{1}{3} \cong \Bi{1}
$$
by \eqref{eq:B1-identification}, and similarly
$$
\faktor{\GammaB{2}{n}}{E^{\mathcal B}_n} \cong \GammaB{2}{3} \cong \Bi{2}
$$
by \eqref{eq:B2-identification}. The geometric statements follow from Proposition~\ref{prop:fibration-principle}.
\end{proof}

\begin{remark}\label{rem:B1B2-covering} 
For every $n\geq3$, $\GammaB{1}{n} \triangleleft \GammaB{2}{n}$, and $\faktor{\GammaB{2}{n}}{\GammaB{1}{n}} \cong\mathbb Z_2$.  Indeed,
$$
\GammaB{2}{n} = \overline P_n\GammaB{1}{n}. 
$$
Moreover, $\overline P_n$ normalizes $\GammaB{1}{n}$. To see this, it suffices to consider conjugation of $\varepsilon$ by elements of $\overline P_n$. The only non-trivial differences produced by the transposition $(1,2)$ are generated by $A_{1,3}A_{2,3}^{-1}=d$ and, for $j\geq4$, by $A_{1,j}A_{2,j}^{-1}\in E^{\mathcal B}_n$. Thus, for every $p\in\overline P_n$, we have $p\varepsilon p^{-1}\varepsilon^{-1}=p((1,2)\cdot p)^{-1}\in\Lambda_{B_1,n}$, proving normality. Consequently,
$$
\faktor{\GammaB{2}{n}}{\GammaB{1}{n}} \cong \faktor{\overline P_n}{\Lambda_{B_1,n}} \cong \mathbb Z_2.
$$

Thus there is a regular double Riemannian covering $X_{\GammaB{1}{n}} \longrightarrow X_{\GammaB{2}{n}}$. For $n=3$ this gives the regular double covering $\XBi{1}\longrightarrow\XBi{2}$. Moreover, for $n\geq4$ one has $\GammaB{2}{n}=\Gammai{2}{n}$ as subgroups of $\overline B_n$. However, the two constructions use distinct complete normal subgroups, namely $\KB{2}{n} = E^{\mathcal B}_n$ and $\Ki{2}{n}=K_{2,4}\oplus E_{2,n}$. Consequently, the same higher-dimensional Bieberbach group admits geometric torus fibrations over both $\XBi{2}$ and $\XGi{2}$.
\end{remark}

\subsection{The groups $\Bi{3}$ and $\Bi{4}$: holonomy $\mathbb Z_2^2$}\label{subsec:B34}

We now consider a Klein four subgroup different from that used in Section~\ref{sec:G6}. Let
$$
a=(1,2), \qquad e=(3,4), \qquad V^-=\langle a,e\rangle\leq S_4.
$$
Then $V^-\cong\mathbb Z_2^2$. Note that $V^-$ contains odd permutations, whereas the Klein four group used for $\Gi{6}$ consists entirely of even permutations. By Proposition~\ref{prop:2groups}, the full preimage
$$
\widetilde V^-=\overline\sigma^{-1}(V^-)
$$
is a Bieberbach group of dimension $6$.

We first construct a $V^-$-invariant full-rank sublattice of $\overline P_4$. Consider the rank-$4$ block $W= \langle A_{1,3},A_{1,4},A_{2,3},A_{2,4} \rangle$. Define
$$
\begin{aligned}
u&=A_{1,3}A_{1,4}A_{2,3}^{-1}A_{2,4}^{-1},\\
v&=A_{1,3}A_{1,4}^{-1}A_{2,3}^{-1}A_{2,4},\\
w_0&=A_{1,3}A_{1,4}A_{2,3}A_{2,4},\\
w_e&=A_{1,3}A_{1,4}^{-1}A_{2,3}A_{2,4}^{-1}.
\end{aligned}
$$
With respect to the ordered basis $(A_{1,3},A_{1,4},A_{2,3},A_{2,4})$ of $W$, the columns of the matrix
$$
H= 
\begin{pmatrix}
 1& 1&1& 1\\
 1&-1&1&-1\\
-1&-1&1& 1\\
-1& 1&1&-1
\end{pmatrix}
$$
are the exponent vectors of $(u,v,w_0,w_e)$. Since $\det(H)=16$, the subgroup $\langle u,v,w_0,w_e\rangle$ is a sublattice of index $16$ in $W$.

Define
\begin{equation}
\Lambda^-_4 = \langle A_{1,2},A_{3,4},u,v,w_0,w_e \rangle. 
\label{eq:B34-lattice}
\end{equation}
Since $\overline P_4 = \langle A_{1,2},A_{3,4}\rangle\oplus W$, we obtain $[\overline P_4:\Lambda^-_4]=16$. The lattice $\Lambda^-_4$ is $V^-$-invariant. Indeed, the action of $a=(1,2)$ is
\begin{equation}
\begin{aligned}
u&\longmapsto u^{-1},
&
v&\longmapsto v^{-1},
&
w_0&\longmapsto w_0,
&
w_e&\longmapsto w_e,
\end{aligned}
\label{eq:B34-action-a}
\end{equation}
whereas the action of $e=(3,4)$ is
\begin{equation}
\begin{aligned}
u&\longmapsto u,
&
v&\longmapsto v^{-1},
&
w_0&\longmapsto w_0,
&
w_e&\longmapsto w_e^{-1}.
\end{aligned}
\label{eq:B34-action-e}
\end{equation}
Both $A_{1,2}$ and $A_{3,4}$ are fixed by $V^-$.

We now choose lifts of the two generators of $V^-$. Put $\alpha=\sigma_1$, $\varepsilon_3=A_{1,3}A_{1,4}\sigma_3$, $\varepsilon_4=A_{1,3}^2\sigma_3$. Then
$$
\overline\sigma(\alpha)=a, \qquad \overline\sigma(\varepsilon_3) = \overline\sigma(\varepsilon_4)=e. 
$$
Since $\alpha=\sigma_1$, we have 
\begin{equation}
\alpha^2=A_{1,2}.
\label{eq:B34-alpha-square}
\end{equation}

The two lifts of $e$ have the same square. Indeed, using \eqref{eq:pair-action},
$$
\varepsilon_3^2 = (A_{1,3}A_{1,4}) \bigl(e\cdot(A_{1,3}A_{1,4})\bigr)A_{3,4} = A_{1,3}^2A_{1,4}^2A_{3,4}, 
$$
whereas
$$
\varepsilon_4^2 = A_{1,3}^2 \bigl(e\cdot A_{1,3}^2\bigr)A_{3,4} = A_{1,3}^2A_{1,4}^2A_{3,4}. 
$$
Since $uw_0=A_{1,3}^2A_{1,4}^2$, we obtain
\begin{equation}
\varepsilon_3^2 = \varepsilon_4^2 = uw_0A_{3,4}. 
\label{eq:B34-epsilon-square}
\end{equation}

The difference between the two extension classes appears in the mixed relation between the two holonomy lifts. Let $\lambda_3=A_{1,3}A_{1,4}$ and $\lambda_4=A_{1,3}^2$, so that $\varepsilon_j=\lambda_j\sigma_3$ for $j=3,4$. Since $\sigma_1$ and $\sigma_3$ commute, one has 
$$
\varepsilon_j\alpha\varepsilon_j^{-1} = \lambda_j \bigl(a\cdot\lambda_j^{-1}\bigr) \alpha. 
$$
For $j=3$, $\lambda_3 \bigl(a\cdot\lambda_3^{-1}\bigr) = A_{1,3}A_{1,4}A_{2,3}^{-1}A_{2,4}^{-1} = u$, and hence
\begin{equation}
\varepsilon_3\alpha\varepsilon_3^{-1} = u\alpha. 
\label{eq:B3-extension-relation}
\end{equation}
For $j=4$, $\lambda_4 \bigl(a\cdot\lambda_4^{-1}\bigr) = A_{1,3}^2A_{2,3}^{-2} = uv$, and therefore
\begin{equation}
\varepsilon_4\alpha\varepsilon_4^{-1} = uv\alpha. 
\label{eq:B4-extension-relation}
\end{equation}

Define
$$
\GammaB{3}{4} = \langle\Lambda^-_4,\alpha,\varepsilon_3\rangle, \qquad \GammaB{4}{4} = \langle\Lambda^-_4,\alpha,\varepsilon_4\rangle. 
$$

\begin{proposition}
\label{prop:B34-bieb}
The groups $\GammaB{3}{4}$ and $\GammaB{4}{4}$ are Bieberbach groups of dimension $6$, with translation lattice $\Lambda^-_4$ and holonomy group $V^-\cong\mathbb Z_2^2$.
\end{proposition}

\begin{proof}
By \eqref{eq:B34-action-a}--\eqref{eq:B34-action-e}, the lattice $\Lambda^-_4$ is normalized by $\alpha$, $\varepsilon_3$ and $\varepsilon_4$. Equations \eqref{eq:B34-alpha-square} and \eqref{eq:B34-epsilon-square} show that the squares of the chosen lifts belong to $\Lambda^-_4$, while \eqref{eq:B3-extension-relation} and \eqref{eq:B4-extension-relation} show that their classes commute modulo $\Lambda^-_4$.

Hence, for $j=3,4$, the quotient $\faktor{\GammaB{j}{4}}{\Lambda^-_4}$ is generated by two commuting involutions and therefore has order at most $4$. On the other hand, $\overline\sigma$ induces a surjection $\faktor{\GammaB{j}{4}}{\Lambda^-_4} \longrightarrow V^-$, and $|V^-|=4$. Thus $\faktor{\GammaB{j}{4}}{\Lambda^-_4} \cong V^-$. It follows in particular that
$$
\GammaB{j}{4}\cap\overline P_4 = \Lambda^-_4, 
$$
so $\Lambda^-_4$ is precisely the translation subgroup of $\GammaB{j}{4}$.

By Proposition~\ref{prop:2groups}, $\widetilde V^-=\overline\sigma^{-1}(V^-)$ is a Bieberbach group. Each $\GammaB{j}{4}$ is a subgroup of $\widetilde V^-$, hence is discrete and torsion-free. Since $\Lambda^-_4$ has full rank in $\overline P_4$, the quotient is compact. Therefore $\GammaB{j}{4}$ is a Bieberbach group of dimension $6$ with the stated holonomy.
\end{proof}

We now isolate the rank-$3$ quotient that produces the two non-orientable Bieberbach groups. Put $K^-_4 = \langle A_{3,4},w_0,w_e \rangle$. The three displayed generators are elements of the basis \eqref{eq:B34-lattice} of $\Lambda^-_4$. Hence $K^-_4$ is primitive of rank $3$ and $\faktor{\Lambda^-_4}{K^-_4} \cong \mathbb Z^3$. Moreover, \eqref{eq:B34-action-a}--\eqref{eq:B34-action-e} show that $K^-_4$ is $V^-$-invariant, and therefore it is normal in both $\GammaB{3}{4}$ and $\GammaB{4}{4}$.

In the quotient lattice, put
$$
t_1=\overline{A}_{1,2}, \qquad t_2=\overline u, \qquad t_3=\overline v. 
$$
Relative to this basis, the two holonomy generators act by
\begin{equation}
\rho(a) = \operatorname{diag}(1,-1,-1), \qquad \rho(e) = \operatorname{diag}(1,1,-1).
\label{eq:B34-holonomy-matrices}
\end{equation}
Thus this is precisely the integral holonomy representation shared by the Bieberbach groups $\Bi{3}$ and $\Bi{4}$.

\begin{theorem}\label{thm:B34-dim4} 
The quotients of the groups in Proposition~\ref{prop:B34-bieb} by $K^-_4$ are 
$$
\faktor{\GammaB{3}{4}}{K^-_4} \cong \Bi{3}, \qquad \faktor{\GammaB{4}{4}}{K^-_4} \cong \Bi{4}.
$$
\end{theorem}

\begin{proof}
First consider $Q_3 = \faktor{\GammaB{3}{4}}{K^-_4}$.  Its translation subgroup is free abelian with basis $t_1,t_2,t_3$. By \eqref{eq:B34-alpha-square} and \eqref{eq:B34-epsilon-square}, $\alpha^2=t_1$, $\varepsilon_3^2=t_2$. The holonomy action \eqref{eq:B34-holonomy-matrices} gives
$$
\alpha t_2\alpha^{-1}=t_2^{-1}, \qquad \alpha t_3\alpha^{-1}=t_3^{-1},
$$
and
$$
\varepsilon_3t_1\varepsilon_3^{-1}=t_1, \qquad \varepsilon_3t_3\varepsilon_3^{-1}=t_3^{-1}.
$$
Finally, \eqref{eq:B3-extension-relation} becomes
$$
\varepsilon_3\alpha\varepsilon_3^{-1} = t_2\alpha. 
$$
Together with commutativity of $t_1,t_2,t_3$, these are precisely the defining relations of the Bieberbach group $\Bi{3}$; see \cite[Theorem~3.5.9]{Wolf} and \cite[List~(A)]{HJKL2002}. Hence $Q_3\cong\Bi{3}$. 

Now consider $Q_4= \faktor{\GammaB{4}{4}}{K^-_4}$. The translation subgroup, holonomy action and square relations are the same as for $Q_3$, but \eqref{eq:B4-extension-relation} becomes
$$
\varepsilon_4\alpha\varepsilon_4^{-1} = t_2t_3\alpha. 
$$
Thus the only change in the extension presentation is the factor $t_3$ in the mixed relation between the two holonomy lifts.  These are precisely the defining relations of the Bieberbach group $\Bi{4}$ in the same classification. Therefore $Q_4\cong\Bi{4}$. 
\end{proof}

We finally stabilize these constructions. For $n\geq4$, set
$$
E^-_n = \langle A_{i,j} \mid 1\leq i<j\leq n,\ j\geq5 \rangle, 
$$
with $E^-_4=\{1\}$. Since $V^-$ acts only on the first four strings, $E^-_n$ is $V^-$-invariant. Define $\Lambda^-_n = \Lambda^-_4\oplus E^-_n$. Then $\Lambda^-_n$ has full rank $N$ in $\overline P_n$, and $[\overline P_n:\Lambda^-_n]=16$. For $j=3,4$, define 
$$
\GammaB{j}{n} = \langle \Lambda^-_n,\alpha,\varepsilon_j \rangle. 
$$
Finally put $K^-_n = K^-_4\oplus E^-_n$. Since $\rank E^-_n=N-6$, one has $\rank K^-_n = 3+(N-6) = N-3$. 

\begin{theorem}\label{thm:B34} 
For every $n\geq4$, the groups $\GammaB{3}{n}$ and $\GammaB{4}{n}$ are Bieberbach subgroups of $\overline B_n$ of dimension $N$ and holonomy $\mathbb Z_2^2$. The subgroup $K^-_n\cong\mathbb Z^{N-3}$ is normal in both groups and 
$$
\faktor{\GammaB{3}{n}}{K^-_n} \cong \Bi{3}, \qquad \faktor{\GammaB{4}{n}}{K^-_n} \cong \Bi{4}. 
$$
Consequently, there are geometric torus bundles $T^{N-3} \longrightarrow X_{\GammaB{3}{n}} \longrightarrow \XBi{3}$ and $T^{N-3} \longrightarrow X_{\GammaB{4}{n}} \longrightarrow \XBi{4}$, with totally geodesic fibers.
\end{theorem}

\begin{proof}
The lattice $\Lambda^-_n$ is invariant under $V^-$, and the relations \eqref{eq:B34-alpha-square}, \eqref{eq:B34-epsilon-square}, \eqref{eq:B3-extension-relation} and \eqref{eq:B4-extension-relation} remain valid in $\overline B_n$.  Exactly as in Proposition~\ref{prop:B34-bieb}, for $j=3,4$ one obtains $\faktor{\GammaB{j}{n}}{\Lambda^-_n} \cong V^-$, and therefore $\Lambda^-_n$ is the translation lattice of $\GammaB{j}{n}$.

Moreover, $\GammaB{j}{n} \leq \overline\sigma^{-1}(V^-)$, and the latter is a Bieberbach group by Proposition~\ref{prop:2groups}.  Hence $\GammaB{j}{n}$ is torsion-free and discrete. Since $\Lambda^-_n$ has full rank in $\overline P_n$, $\GammaB{j}{n}$ is cocompact and therefore a Bieberbach group of dimension $N$.

The subgroup $K^-_n$ is normal by construction and has rank $N-3$. Quotienting by $E^-_n$ reduces the construction to the case $n=4$. Thus Theorem~\ref{thm:B34-dim4} gives
$$
\faktor{\GammaB{3}{n}}{K^-_n} \cong \Bi{3}, \qquad \faktor{\GammaB{4}{n}}{K^-_n} \cong \Bi{4}. 
$$
The geometric statements follow from Proposition~\ref{prop:fibration-principle}.
\end{proof}

\section{The realization theorem and splitting of the torus bundles}\label{sec:realization-splitting}

With the case-by-case constructions now in place, we first assemble them to prove Theorem~\ref{thm:main}, and then determine which of the associated torus-bundle extensions split.

\subsection{Proof of the main theorem}\label{sec:main-proof}

The optimality of the thresholds was established in Remark~\ref{rmk:mainthm}; it remains to assemble the constructions.

\begin{proof}[Proof of Theorem~\ref{thm:main}]
Let $M$ be a compact connected flat $3$-manifold. By Wolf's classification \cite[Theorems~3.5.5 and~3.5.9]{Wolf}, its fundamental group is isomorphic to one of
$$
\Gi{1},\ldots,\Gi{6}, \Bi{1},\ldots,\Bi{4}.
$$

For $\Gi{1}\cong\mathbb Z^3$, the direct-summand construction in the Introduction gives the required realization for every $n\geq3$.

For the remaining orientable cases, the required constructions are provided by Theorems~\ref{thm:G2}, \ref{thm:G3}, \ref{thm:G4}, \ref{thm:G5}, and~\ref{thm:G6}, with thresholds $4,3,4,5$, and~$4$, respectively. For the non-orientable cases, Theorem~\ref{thm:B12} treats both $\Bi{1}$ and $\Bi{2}$ for every $n\geq3$, while Theorem~\ref{thm:B34} treats both $\Bi{3}$ and $\Bi{4}$ for every $n\geq4$.

Thus, in every case, there exist a Bieberbach subgroup $\Gamma_n(M)\leq\overline B_n$ of dimension $N$ and a normal free abelian subgroup $K_n(M)$ of rank $N-3$, contained in the translation lattice of $\Gamma_n(M)$, such that
$$
\faktor{\Gamma_n(M)}{K_n(M)}\cong\pi_1(M).
$$
Since the quotient is a $3$-dimensional Bieberbach group, Proposition~\ref{prop:fibration-principle} implies that $K_n(M)$ is complete and yields the geometric torus bundle
$$
T^{N-3}\longrightarrow X_{\Gamma_n(M)}\longrightarrow X_{\faktor{\Gamma_n(M)}{K_n(M)}}
$$
with totally geodesic fibers. Since $\faktor{\Gamma_n(M)}{K_n(M)}\cong\pi_1(M)$, it follows from \cite[Theorem~3.3.1]{Wolf} (see also \cite[Chapter~I, Theorem~4.1]{Charlap}) that $X_{\faktor{\Gamma_n(M)}{K_n(M)}}$ is affinely equivalent to $M$. This proves the theorem.
\end{proof}

\subsection{Splitting of the torus-bundle extensions}\label{sec:splitting}

Each realization $\big(\Gamma_n(M),K_n(M)\big)$ constructed in the Introduction and in Sections~\ref{sec:cyclic} and~\ref{sec:nonorientable} determines a short exact sequence
\begin{equation}
1 \longrightarrow K_n(M) \longrightarrow \Gamma_n(M) \longrightarrow \pi_1(M) \longrightarrow 1,
\label{eq:bundle-extension}
\end{equation}
where $K_n(M)\cong\mathbb Z^{N-3}$ is the fundamental group of the torus fiber. We say that \eqref{eq:bundle-extension} splits if the quotient homomorphism admits a group-theoretic section. We use two elementary lemmas.

\begin{lemma}\label{lem:splitting-obstructions}
Let $\Gamma$ be a Bieberbach group with translation lattice $\Lambda$ and holonomy group $F$, and let $K\leq\Lambda$ be a normal subgroup of $\Gamma$ such that $\Pi=\faktor{\Gamma}{K}$ is a Bieberbach group whose holonomy group has order $|F|$. Let $L$ be the translation lattice of $\Pi$, and let $p\colon\Gamma\longrightarrow\Pi$ be the quotient map. Then the following hold.
\begin{itemize}
\item[\textup{(a)}] The restriction $p|_\Lambda\colon\Lambda\longrightarrow L$ is surjective with kernel $K$, and $p$ induces an isomorphism from $F$ onto the holonomy group of $\Pi$. Identifying these two groups, $p|_\Lambda$ is a homomorphism of $\mathbb Z[F]$-modules.
\item[\textup{(b)}] If $p$ admits a group-theoretic section, then the exact sequence of $\mathbb Z[F]$-modules
$$
0\longrightarrow K\longrightarrow\Lambda\xrightarrow{\;p|_\Lambda\;} L\longrightarrow0
$$
splits.
\end{itemize}
\end{lemma}

\begin{proof}
(a) The first two assertions are Lemma~\ref{lem:holonomy-quotient} in the case of equal holonomy orders. If $\gamma\in\Gamma$ has holonomy $g$ and $\lambda\in\Lambda$, then $p(\gamma\lambda\gamma^{-1})=p(\gamma)p(\lambda)p(\gamma)^{-1}$, and $p(\gamma)$ has holonomy corresponding to $g$; hence $p|_\Lambda$ is $F$-equivariant.

(b) Let $\tau\colon\Pi\longrightarrow\Gamma$ be a section of $p$. For $\ell\in L$, the holonomy of $\tau(\ell)$ corresponds, under the isomorphism in (a), to the holonomy of $p(\tau(\ell))=\ell$, which is trivial. Hence $\tau(\ell)\in\Lambda$, and $\tau|_L\colon L\longrightarrow\Lambda$ is a homomorphism with $p\circ\tau|_L=\mathrm{id}_L$. Let $g\in F$, and choose $\gamma\in\Pi$ with holonomy corresponding to $g$. By the same argument, $\tau(\gamma)$ has holonomy $g$ in $\Gamma$, and therefore
$$
\tau(g\cdot\ell)=\tau(\gamma\ell\gamma^{-1})=\tau(\gamma)\tau(\ell)\tau(\gamma)^{-1}=g\cdot\tau(\ell).
$$
Thus $\tau|_L$ is an $F$-equivariant section of $p|_\Lambda$, and $\Lambda=K\oplus\tau(L)$ as $\mathbb Z[F]$-modules.
\end{proof}

All the extensions~\eqref{eq:bundle-extension} satisfy the hypotheses of Lemma~\ref{lem:splitting-obstructions}: in each case $K_n(M)$ consists of translations, and $\Gamma_n(M)$ and $\pi_1(M)$ have the same holonomy group, namely the trivial group for $\Gi{1}$, $\mathbb Z_2$ for $\Gi{2},\Bi{1},\Bi{2}$, $\mathbb Z_3$ for $\Gi{3}$, $\mathbb Z_4$ for $\Gi{4}$, $\mathbb Z_6$ for $\Gi{5}$, and $\mathbb Z_2^2$ for $\Gi{6},\Bi{3},\Bi{4}$.

The second lemma reduces each stabilized construction to its seed case.

\begin{lemma}\label{lem:seed-reduction}
Let $3\leq m\leq n$, let $\Gamma_m\leq\overline B_m$, and let $K_m$ be a normal subgroup of $\Gamma_m$. Put $E=E^{(m)}_n$, $\Gamma_n=E\,\Gamma_m$ and $K_n=E\,K_m$. Then $K_n$ is normal in $\Gamma_n$, and if the extension $1\to K_n\to\Gamma_n\to\faktor{\Gamma_n}{K_n}\to1$ splits, so does $1\to K_m\to\Gamma_m\to\faktor{\Gamma_m}{K_m}\to1$.
\end{lemma}

\begin{proof}
By Remark~\ref{rem:stabilization}, $\Gamma_n=E\rtimes\Gamma_m$. Let $\pi\colon\Gamma_n\longrightarrow\Gamma_m$ be the retraction with kernel $E$, given by $\pi(e\gamma)=\gamma$ for $e\in E$ and $\gamma\in\Gamma_m$. Then $K_n=\pi^{-1}(K_m)$, so $K_n$ is normal in $\Gamma_n$, and $\pi$ induces an isomorphism
$$
\overline\pi\colon\faktor{\Gamma_n}{K_n}\longrightarrow\faktor{\Gamma_m}{K_m}
\qquad\text{with}\qquad
p_m\circ\pi=\overline\pi\circ p_n,
$$
where $p_n$ and $p_m$ are the quotient maps. If $\tau$ is a section of $p_n$, then $\pi\circ\tau\circ\overline\pi^{-1}$ is a section of $p_m$, since $p_m\circ\pi\circ\tau\circ\overline\pi^{-1}=\overline\pi\circ p_n\circ\tau\circ\overline\pi^{-1}=\mathrm{id}$.
\end{proof}

Directly from the definitions in Sections~\ref{sec:cyclic} and~\ref{sec:nonorientable}, every non-trivial realization is of the form $\Gamma_n(M)=E^{(m)}_n\,\Gamma_m(M)$ and $K_n(M)=E^{(m)}_n\,K_m(M)$, with the following seed data:
$$
\begin{array}{c|c|c|c|c}
\pi_1(M) & m & E^{(m)}_n & \Gamma_m(M) & K_m(M)\\ \hline
\Gi{2} & 4 & E_{2,n} & \Gammai{2}{4} & K_{2,4}\\
\Gi{3} & 3 & E_{3,n} & \Gammai{3}{3} & \{1\}\\
\Gi{4} & 4 & E_{4,n} & \Gammai{4}{4} & K_{4,4}\\
\Gi{5} & 5 & E_{5,n} & \Gammai{5}{5} & K_{5,5}\\
\Gi{6} & 4 & E^{(4)}_n & \Gammai{6}{4} & K^{(6)}_4\\
\Bi{1},\ \Bi{2} & 3 & E^{\mathcal B}_n & \GammaB{1}{3},\ \GammaB{2}{3} & \{1\}\\
\Bi{3},\ \Bi{4} & 4 & E^{-}_n & \GammaB{3}{4},\ \GammaB{4}{4} & K^-_4
\end{array}
$$
For instance, $\Gammai{2}{n}=\GammaB{2}{n}=\overline P_n\langle\sigma_1\rangle=E^{(m)}_n\,\overline P_m\langle\sigma_1\rangle$ with $m=4$ and $m=3$, respectively, and for $\Gi{6}$ the generating set $\mathcal B_{6,n}$ shows that $\Lambda_{6,n}=\Lambda_{6,4}\oplus E^{(4)}_n$ and $K^{(6)}_n=K^{(6)}_4\oplus E^{(4)}_n$.

\begin{theorem}\label{thm:splitting}
For each compact flat $3$-manifold $M$ and every $n\geq n_M$, the extension~\eqref{eq:bundle-extension} associated with the realization $\big(\Gamma_n(M),K_n(M)\big)$ constructed above has the following splitting behavior, according to the Bieberbach group $\pi_1(M)$:
$$
\begin{array}{c|l}
\text{splits}
& \Gi{1},\ \Gi{3},\ \Bi{1},\ \Bi{2} \\[1mm]
\text{does not split} & \Gi{2},\ \Gi{4},\ \Gi{5},\ \Gi{6},\ \Bi{3},\ \Bi{4}.
\end{array}
$$
\end{theorem}

\begin{proof}
Throughout, translation lattices are written additively when we use module structures, and multiplicatively when we work inside the groups.

\medskip

\noindent
\textit{The split cases.}
For $\Gi{1}$, the total group is $\overline P_n\cong\mathbb Z^N$, and $K$ was chosen to be a direct summand with quotient $\mathbb Z^3$. Hence the extension splits.

For $\Gi{3}$, $\Bi{1}$ and $\Bi{2}$, the seed data above give $m=3$ and $K_3(M)=\{1\}$, so that $K_n(M)=E^{(3)}_n$. By Remark~\ref{rem:stabilization},
$$
\Gammai{3}{n}=E^{(3)}_n\rtimes\Gammai{3}{3},\qquad \GammaB{1}{n}=E^{(3)}_n\rtimes\GammaB{1}{3},\qquad \GammaB{2}{n}=E^{(3)}_n\rtimes\GammaB{2}{3}.
$$
In each case the complement $\Gamma_3(M)$ maps isomorphically onto $\faktor{\Gamma_n(M)}{K_n(M)}$, and the inverse of this isomorphism is a section. Hence these three extensions split.

\medskip

\noindent
For the six remaining cases, Lemma~\ref{lem:seed-reduction} shows that it suffices to prove that the seed extension $1\to K_m(M)\to\Gamma_m(M)\to\pi_1(M)\to1$ does not split.

\medskip

\noindent
\textit{The case $\Gi{2}$.}
Here $m=4$, $\Gamma=\Gammai{2}{4}=\overline P_4\langle\sigma_1\rangle$, $\Lambda=\overline P_4$, and $F\cong\mathbb Z_2$ is generated by the class of $\sigma_1$. By~\eqref{eq:pair-action}, $\sigma_1$ acts on $\overline P_4$ through the transposition $(1,2)$: it fixes $A_{1,2}$ and $A_{3,4}$, and interchanges $A_{1,3}$ with $A_{2,3}$ and $A_{1,4}$ with $A_{2,4}$. Let
$$
\kappa\colon\overline P_4\longrightarrow L=\faktor{\overline P_4}{K_{2,4}}
$$
be the quotient map, and put $u=\kappa(A_{1,3})$ and $v=\kappa(A_{1,4})$. Since $A_{1,3}+A_{2,3}$, $A_{1,4}+A_{2,4}$ and $A_{3,4}$ form part of a basis of $\overline P_4$ together with $A_{1,2}$, $A_{1,3}$, $A_{1,4}$, the lattice $L$ is free abelian with basis $\kappa(A_{1,2}),u,v$, and $\kappa(A_{2,3})=-u$, $\kappa(A_{2,4})=-v$. In particular, $\sigma_1\cdot u=-u$.

Suppose that the extension splits. By Lemma~\ref{lem:splitting-obstructions}, there is a $\mathbb Z[\mathbb Z_2]$-module section $\iota\colon L\longrightarrow\overline P_4$ of $\kappa$. Put $U=\iota(u)=\sum_{i<j}c_{i,j}A_{i,j}$. Then
$$
\sigma_1\cdot U=\iota(\sigma_1\cdot u)=-U,
$$
which forces $c_{1,2}=c_{3,4}=0$, $c_{2,3}=-c_{1,3}$ and $c_{2,4}=-c_{1,4}$. Hence
$$
U=c_{1,3}(A_{1,3}-A_{2,3})+c_{1,4}(A_{1,4}-A_{2,4}),
\qquad
\kappa(U)=2c_{1,3}\,u+2c_{1,4}\,v\in 2L.
$$
This contradicts $\kappa(U)=u$, since $u$ is a basis element of $L$. Hence the $\Gi{2}$ extension does not split.

\medskip

\noindent
\textit{The case $\Gi{6}$.}
Here $m=4$, and we write $K=K^{(6)}_4=\langle s_1,s_2,s_3\rangle$ and $p\colon\Gammai{6}{4}\longrightarrow\faktor{\Gammai{6}{4}}{K}\cong\Gi{6}$ for the quotient map. Suppose that $\tau$ is a section of $p$.

Since $\Lambda_{6,4}$ is abelian, conjugation by any element of $\beta K$ acts on $\Lambda_{6,4}$ as conjugation by $\beta$, which by~\eqref{eq:G6-3action} inverts $d_1$ and fixes $K$ pointwise; similarly, conjugation by any element of $\alpha K$ fixes $K$ pointwise. Now $\tau(p(d_1))\in p^{-1}(p(d_1))=d_1K$, say $\tau(p(d_1))=d_1k$ with $k\in K$, and $\tau(p(\beta))\in\beta K$. Applying $\tau$ to the relation $p(\beta)p(d_1)p(\beta)^{-1}=p(d_1)^{-1}$ gives
$$
d_1^{-1}k=(d_1k)^{-1}=d_1^{-1}k^{-1}.
$$
Hence $k^2=1$, and since $K$ is torsion-free, $\tau(p(d_1))=d_1$.

Next, $\tau(p(\alpha))=h\alpha$ for some $h\in K$. Since $\alpha$ fixes $K$ pointwise, \eqref{eq:G6-alpha2} gives
$$
\tau(p(\alpha))^2=(h\alpha)^2=h^2\alpha^2=h^2d_1s_1s_2s_3^2.
$$
On the other hand, $p(\alpha)^2=p(d_1s_1s_2s_3^2)=p(d_1)$, so $\tau(p(\alpha))^2=\tau(p(d_1))=d_1$. Therefore
$$
h^2=(s_1s_2s_3^2)^{-1},
$$
which is impossible in the free abelian group $K$, since the exponents of $s_1$ and $s_2$ on the right-hand side are odd. Thus the $\Gi{6}$ extension does not split.

\medskip

\noindent
\textit{The case $\Gi{4}$.}
Here $m=4$, $\Lambda=\Lambda_{4,4}$ with basis $x_0,x_1,x_2,x_3,s,d$, and $F\cong\mathbb Z_4$ is generated by the class of $\theta$, which acts by $x_r\mapsto x_{r+1}$ $(r\in\mathbb Z/4\mathbb Z)$, $s\mapsto s$ and $d\mapsto -d$. Let
$$
\kappa\colon\Lambda_{4,4}\longrightarrow L=\faktor{\Lambda_{4,4}}{K_{4,4}}
$$
be the quotient map. By~\eqref{eq:K44}, $L$ is free abelian with basis $u=\kappa(x_0)$, $v=\kappa(x_1)$, $z=\kappa(s)$, and
$$
\kappa(x_2)=-u,\qquad \kappa(x_3)=-v,\qquad \kappa(d)=0.
$$
In particular, $\theta^2\cdot u=-u$; compare~\eqref{eq:G4-action-quot}.

Suppose that the extension splits. By Lemma~\ref{lem:splitting-obstructions}, there is a $\mathbb Z[\mathbb Z_4]$-module section $\iota\colon L\longrightarrow\Lambda_{4,4}$ of $\kappa$. Put
$$
U=\iota(u)=a_0x_0+a_1x_1+a_2x_2+a_3x_3+bs+cd.
$$
Equivariance gives $\theta^2\cdot U=\iota(\theta^2\cdot u)=-U$. Since $\theta^2$ interchanges $x_0$ with $x_2$ and $x_1$ with $x_3$, and fixes $s$ and $d$, we obtain
$$
a_2=-a_0,\qquad a_3=-a_1,\qquad b=c=0.
$$
Consequently,
$$
\kappa(U)=(a_0-a_2)u+(a_1-a_3)v=2a_0u+2a_1v\in 2L,
$$
contradicting $\kappa(U)=u$. Hence the $\Gi{4}$ extension does not split.

\medskip

\noindent
\textit{The case $\Gi{5}$.}
Here $m=5$, $\Lambda=\mathcal L_{5,5}$, and $F\cong\mathbb Z_6$ is generated by the class of $\theta$. With the notation of Section~\ref{subsec:G5}, there is a decomposition into $F$-invariant sublattices
$$
\mathcal L_{5,5}=\big(L_0\oplus\langle t\rangle\big)\oplus W.
$$
Indeed, $\theta$ fixes $s$ and $t$, permutes $A_{1,2}^3=s^3p^{-1}q^{-1}$, $A_{1,3}^3=p$ and $A_{2,3}^3=q$, and permutes $w_0,\ldots,w_5$ cyclically. Moreover,
$$
K_{5,5}=K'\oplus K_W,\qquad K'=\langle ps^{-1},\,qs^{-1},\,st\rangle\leq L_0\oplus\langle t\rangle,
$$
and $\faktor{(L_0\oplus\langle t\rangle)}{K'}$ is infinite cyclic, generated by the class of $s$, with trivial $F$-action. Identify $W$ with $\mathbb Z[x]/(x^6-1)$ via $w_j\mapsto x^j$, so that $\theta$ acts by multiplication by $x$, and put $C=\mathbb Z[x]/(\Phi_6(x))$, where $\Phi_6(x)=x^2-x+1$. Then the quotient map $\kappa\colon\mathcal L_{5,5}\longrightarrow L=\faktor{\mathcal L_{5,5}}{K_{5,5}}$ induces an $F$-equivariant decomposition
$$
L=\langle\overline s\rangle\oplus C,
$$
and the restriction of $\kappa$ to $W$ is the map $q_6$ of~\eqref{eq:G5-cyclotomic-map}, followed by the inclusion of $C$. The module $C$ is free abelian with basis $\{1,x\}$.

Suppose that the extension splits. By Lemma~\ref{lem:splitting-obstructions}, there is a $\mathbb Z[\mathbb Z_6]$-module section $\iota\colon L\longrightarrow\mathcal L_{5,5}$ of $\kappa$. Put $U=\iota(1)$, with $1\in C$, and write $U=U'+Y$ with $U'\in L_0\oplus\langle t\rangle$ and $Y\in W$. Since $\kappa(U')\in\langle\overline s\rangle$, comparing components in $\kappa(U)=1$ gives
$$
q_6(Y)=1.
$$
Since $\Phi_6(x)$ annihilates $C$, we have $\Phi_6(\theta)\cdot U=\iota\big(\Phi_6(\theta)\cdot1\big)=0$, and the $F$-invariance of the decomposition gives $\Phi_6(x)\,Y=0$ in $\mathbb Z[x]/(x^6-1)$.

Let $f\in\mathbb Z[x]$ represent $Y$, and put
$$
Q(x)=\frac{x^6-1}{\Phi_6(x)}=x^4+x^3-x-1.
$$
Then $\Phi_6f=(x^6-1)g=\Phi_6Qg$ for some $g\in\mathbb Z[x]$, and since $\mathbb Z[x]$ is an integral domain, $f=Qg$. Modulo $\Phi_6(x)$ one has $x^2=x-1$, $x^3=-1$ and $x^4=-x$, hence
$$
Q(x)\equiv-2(x+1)\pmod{\Phi_6(x)}.
$$
It follows that $q_6(Y)=-2(x+1)g\in 2C$, which contradicts $q_6(Y)=1$, since $1$ is a basis element of $C$. Hence the $\Gi{5}$ extension does not split.

\medskip

\noindent
\textit{The cases $\Bi{3}$ and $\Bi{4}$.}
Here $m=4$. Fix $j\in\{3,4\}$, let
$$
\pi_j\colon\GammaB{j}{4}\longrightarrow\faktor{\GammaB{j}{4}}{K^-_4}\cong\Bi{j}
$$
be the quotient map, and put $\overline\alpha=\pi_j(\alpha)$, $\overline\varepsilon_j=\pi_j(\varepsilon_j)$ and $t_2=\pi_j(u)$. Suppose that $\tau$ is a section of $\pi_j$.

Since $\Lambda^-_4$ is abelian, conjugation by any element of $\alpha K^-_4$ acts on $\Lambda^-_4$ as conjugation by $\alpha$, which by~\eqref{eq:B34-action-a} inverts $u$ and fixes $K^-_4=\langle A_{3,4},w_0,w_e\rangle$ pointwise. Now $\tau(t_2)=uk$ for some $k\in K^-_4$, and $\tau(\overline\alpha)\in\alpha K^-_4$. Applying $\tau$ to the relation $\overline\alpha t_2\overline\alpha^{-1}=t_2^{-1}$ gives
$$
u^{-1}k=(uk)^{-1}=u^{-1}k^{-1},
$$
so $k^2=1$. As $K^-_4$ is torsion-free, $\tau(t_2)=u$.

Next, $\tau(\overline\varepsilon_j)=h\varepsilon_j$ for some $h=A_{3,4}^p w_0^q w_e^r\in K^-_4$. Since $\varepsilon_j=\lambda_j\sigma_3$ with $\lambda_j\in\overline P_4$, conjugation by $\varepsilon_j$ acts on $\Lambda^-_4$ through $e=(3,4)$, and~\eqref{eq:B34-action-e} gives $e\cdot h=A_{3,4}^p w_0^q w_e^{-r}$. Using~\eqref{eq:B34-epsilon-square}, we obtain, independently of $j$,
$$
(h\varepsilon_j)^2=h\,(e\cdot h)\,\varepsilon_j^2=A_{3,4}^{2p}w_0^{2q}\cdot uw_0A_{3,4}=uA_{3,4}^{2p+1}w_0^{2q+1}.
$$
On the other hand, $\overline\varepsilon_j^{\,2}=\pi_j(uw_0A_{3,4})=t_2$, so $(h\varepsilon_j)^2=\tau(t_2)=u$. Therefore
$$
A_{3,4}^{2p+1}w_0^{2q+1}=1,
$$
which is impossible, since $A_{3,4}$ and $w_0$ belong to the basis~\eqref{eq:B34-lattice} of $\Lambda^-_4$. Hence neither the $\Bi{3}$ extension nor the $\Bi{4}$ extension splits.

This completes the proof.
\end{proof}

\begin{remark}\label{rem:splitting-comparison}
By \cite[Theorem~17]{RT2010}, since $K_n(M)$ is a complete, torsion-free normal subgroup of $\Gamma_n(M)$, the extension~\eqref{eq:bundle-extension} splits if and only if the corresponding geometric torus-bundle projection
$$
X_{\Gamma_n(M)}\longrightarrow X_{\faktor{\Gamma_n(M)}{K_n(M)}}
$$
admits an affine section. Thus Theorem~\ref{thm:splitting} also provides a geometric classification: the torus bundles constructed over $\XGi{1}$, $\XGi{3}$, $\XBi{1}$ and $\XBi{2}$ admit affine sections, whereas those over the remaining six compact flat $3$-manifolds do not.

The splitting behavior depends on the chosen complete normal subgroup, and not only on the total group: for $n\geq4$ one has $\GammaB{2}{n}=\Gammai{2}{n}$, but the extension over $\XBi{2}$ splits whereas the one over $\XGi{2}$ does not. We do not know whether, for the six non-split cases, a different realization of $\pi_1(M)$ in some $\overline B_n$ could give a split extension.
\end{remark}


\begin{thebibliography}{99}

\bibitem{BGH} A.~Bauval, D.~L.~Gon\c{c}alves and C.~Hayat, The Borsuk--Ulam theorem for the Seifert manifolds having flat geometry, arXiv:1807.00159.

\bibitem{Charlap} L.~S.~Charlap, \newblock \emph{Bieberbach Groups and Flat Manifolds}, \newblock Universitext, Springer-Verlag, New York, 1986.

\bibitem{GGO2017} D.~L.~Gon\c{c}alves, J.~Guaschi and O.~Ocampo, \newblock A quotient of the Artin braid groups related to crystallographic groups, \newblock \emph{J. Algebra} \textbf{474} (2017), 393--423.

\bibitem{HJKL2002} K.~Y.~Ha, J.~H.~Jo, S.~W.~Kim and J.~B.~Lee, \newblock Classification of free actions of finite groups on the $3$-torus, \newblock \emph{Topology Appl.} \textbf{121} (2002), 469--507.

\bibitem{KLL2005} S.~W.~Kim, J.~B.~Lee and K.~B.~Lee, Averaging formula for Nielsen numbers, \emph{Nagoya Math. J.} \textbf{178} (2005), 37--53.

\bibitem{Ocampo2021} O.~Ocampo, \newblock Bieberbach groups and flat manifolds with finite abelian holonomy from Artin braid groups, \newblock \emph{J. Knot Theory Ramifications} \textbf{30} (2021), no.~9, 2150069.

\bibitem{OR2019} O.~Ocampo and J.~G.~Rodr\'iguez-Nieto, \newblock On Bieberbach subgroups of $B_n/[P_n,P_n]$ and flat manifolds with cyclic holonomy $\mathbb Z_{2^d}$, \newblock \emph{Topology Appl.} \textbf{265} (2019), 106827.

\bibitem{RT2010} J.~G.~Ratcliffe and S.~T.~Tschantz, \newblock Fibered orbifolds and crystallographic groups, \newblock \emph{Algebr. Geom. Topol.} \textbf{10} (2010), 1627--1664.

\bibitem{Szczepanski} A.~Szczepa\'nski, \newblock \emph{Geometry of Crystallographic Groups}, \newblock Algebra and Discrete Mathematics, vol.~4, World Scientific, Singapore, 2012.

\bibitem{Wolf} J.~A.~Wolf, \newblock \emph{Spaces of Constant Curvature}, \newblock 6th ed., AMS Chelsea Publishing, Providence, RI, 2011. 

\end{thebibliography}
\end{document}